\documentclass[12pt]{amsart}

\usepackage{amssymb,amsfonts}
\usepackage{amsmath,amsthm}
\usepackage{mathtools}
\usepackage{hyperref}
\usepackage{xcolor}
\usepackage{tensor}
\usepackage{tikz-cd}

\usepackage{tikz}
\usetikzlibrary{shapes,snakes,arrows, chains, positioning, shapes.geometric, shapes.symbols,calc}
\usetikzlibrary{arrows.meta}

\def\centerarc[#1](#2)(#3:#4:#5)
{ \draw[#1] ($(#2)+({#5*cos(#3)},{#5*sin(#3)})$) arc (#3:#4:#5); }

\usepackage[margin=3.2cm]{geometry}
\usepackage{imakeidx}
\makeindex

\theoremstyle{plain}
\newtheorem{theorem}{Theorem}

\newtheorem{conjecture}{Conjecture}
\newtheorem{conjpic}{Conjectural picture}
\newtheorem{observation}{Observation}

\theoremstyle{definition}

\newtheorem{remark}{Remark}

\newtheorem*{problem}{Problem}
\newtheorem*{genproblem}{General problem}

\newcommand{\C}{\mathbb{C}}
\newcommand{\R}{\mathbb{R}}

\title[On multiscale aspects in Algebraic Geometry]{On multiscale aspects of K\"ahler-Einstein metrics and Algebraic Geometry}
\author{Cristiano Spotti}

\address{Department of Mathematics, Aarhus University,  Ny Munkegade 118 8000 Aarhus  (Denmark), \emph{c.spotti@math.au.dk}}

\begin{document}

\maketitle

\begin{abstract}
    We discuss how metric limits and rescalings of K\"ahler-Einstein metrics connect with Algebraic Geometry, mostly in relation to the study of moduli spaces of varieties, and singularities. Along the way, we describe some elementary examples, review some recent results, and propose some tentative conjectural pictures.
\end{abstract}

\tableofcontents

\section{Introduction}

K\"ahler–Einstein metrics play a central role in exploring the interplay between Geometric Analysis and Algebraic Geometry: such metrics are a fundamental example of solutions of a non-linear PDE (the Einstein equation) appearing in the study of moduli spaces in Algebraic Geometry, via the so-called Yau-Tian-Donaldson conjecture  (e.g., \cite{chen_donaldson_sun2015_III}). This conjecture is a Hitchin--Kobayashi–type correspondence, which indeed provides a canonical way for ``geometrizing''  algebraic varieties.\\

In this survey, I will focus in particular on the relation between the study of singularities of K\"ahler–Einstein (KE) metrics and Algebraic Geometry, with special emphasis on the role played by the \textit{different scalings} 
$\lambda^2 g_t$ of these canonical metrics. From an analytic point of view, one can interpret such scaling operations as a way to get geometric understanding of the regularity of solutions of these non-linear PDEs, via  ``blowing-up analysis''. However, corresponding notions of  ``canonical scalings'' in Algebraic Geometry remain still quite mysterious, although certain aspects of them have begun to emerge in recent literature, as I will go to discuss in some situations. \\

Thus, the main problem I want to address is the following:

\begin{genproblem}
Given KE metrics $g_t$, understand scalings $\lambda^2 g_t$, linking with Algebraic Geometry.   
\end{genproblem}

I will discuss aspects of such a problem through the study of certain \textit{examples} of KE metrics on polarized algebraic varieties, with two main interrelated objectives in mind:\\

\begin{enumerate}
\item Obtain \textit{computability}, using concrete algebro-geometric tools, of limits and blow-up analysis (regularity) for such KE metrics;\\
\item Give new point of views in the algebro-geometric study of \textit{singularities} and \textit{moduli compactifications}.\\
\end{enumerate}

The first goal is more analytical/differential geometric in nature. Indeed, for general Einstein metrics little is known about a concrete analysis of singularities' formations, beside the deep, but very abstract, regularity theory provided by works of Cheeger and Colding (e.g., \cite{cheeger_degen_riemannian_2001}). However, in the K\"ahler situation we can try to relate such a-priori analytic regularity with the underneath algebraic structures, thanks to the foundational works of Donaldson and Sun \cite{donaldson_sun_gh_limits_2014, donaldson_sun_gh_limits_II_2017}. This ``algebraisation'' can lead to situations where one can \textit{concretely} compute the abstract analytic limits, hence providing a full understanding of solutions for such geometric PDE.   A first example of this philosophy can be given by the study of KE metrics on del Pezzo surfaces (\cite{odaka2016compact}), where one can indeed use Algebraic Geometry to give a completely explicit \textit{classifications} of such metrics together with their Gromov-Hausdorff degenerations  (in particular, completing works by Tian \cite{tian1990calabi} and Mabuchi-Mukai \cite{MabuchiMukai1993}), by reducing  their study to concrete and computable finite dimensional GIT type quotients, via the so-called \textit{moduli continuity method} (e.g., \cite{spotti2019kaehler}).  Combined with results by LeBrun \cite{LeBrun2015}, this gives the first known examples of Gromov-Hausdorff compactifications of  \textit{connected components} of moduli spaces of positive Einstein metrics on the underneath differentiable manifolds. 

Conversely, the second goal is more algebro-geometric. The canonicity of the KE metric structure should point out to peculiar algebro-geometric features for algebraic varieties and singularities that could not be easily visible via already established algebro-geometric frameworks. An example of this is given by the \textit{general} construction of moduli spaces of (K-polystable) Fano varieties, which was indeed made possible thanks to several inputs and ideas coming from KE geometry (in an essential way in the first analytic constructions in the smoothable case \cite{SpottiSunYao2016, Odaka2015CompactModuli,LiWangXu2019}, but still conceptually behind the more recent purely algebro-geometric construction by Xu and collaborators, e.g., \cite{LiuXuZhuang2022}).\\

The KE Fano examples mentioned above are \textit{non-collapsing at fixed volume \textit{scale}} (or fixed diameter, by positivity of the Ricci tensor). Thus, these KE Fano K-moduli spaces are fully describing the geometry of such metrics at such scale. However, it is still possible to ask what happens when one \textit{further rescales} the metrics near singularities (formations). \\

The next interesting scale detects the infinitesimal geometry near the singularities of the degenerate limit $X_{\infty}$ added in the Gromov-Hausdorff compactification, via the notion of \textit{metric tangent cone} \cite{donaldson_sun_gh_limits_II_2017}. Remarkably, such infinitesimal geometry can also be computed algebraically \cite{LiWangXu2021TangentCones}  via the Li's normalized volume $\hat{vol}$, a new local invariant for klt (Kawamata-log-terminal) singularities \cite{Li2018NormalizedVolumes}, whose origin can be traced back to the physics works on AdS-CFT correspondence by Martelli-Sparks-Yau \cite{martelli2006geometric}. 

However, this is not the end of the story from an analytic viewpoint: during the process of formation of singularities one can rescale even further, producing asymptotically conical (AC) Calabi-Yau varieties, still possibly singular. Iterating such higher rescalings one can get a so-called \textit{metric bubble tree} made of affine AC Calabi-Yau varieties, canonically associated to a singularity formation and fully describing the degeneration process. This provides a complete regularity picture for the degenerating family of KE metrics and hence it can be seen as an analytic way to get a canonical resolution of singularities for such geometric shapes.  Thus, a natural question arises: can these trees be computed algebro-geometrically as well? Motivated by some examples, I will discuss some conjectural picture \ref{CP1}, as emerged in a series of related works by  Sun \cite{SSBubbling}, by myself in a joint work with de Borbon \cite{de2024some} and by  Odaka \cite{Odaka2024Bubbling}, aiming to address this aspect.\\

For KE metrics with non positive constant scalar curvature (so for Calabi-Yau or positive canonical bundle varieties) \textit{collapsing} can happen, making the multiscale study of the degenerating geometry more challenging. Local collapsing happens for negative KE metrics. In such case, \textit{generic} limits at constant volume scale are essentially understood in relation to Algebraic Geometry via the KSBA moduli compactification, e.g., \cite{Kollar2013Dollar, SongSturmWang2020}, generalizing the classical relation between the Deligne-Mumford compactification and degenerations of hyperbolic metrics on Riemann surfaces. However, a full multiscale metric understanding near the cuspidal regions remains in general more mysterious (but see, e.g., \cite{Kobayashi1985, ZhangKE, fu2025continuous, DatarFuSong2023}, and some conjectural picture for special situations I will propose at the end of this survey \ref{CP3}). 

For polarized Calabi-Yau varieties, instead, collapsing happens \textit{everywhere} on the manifolds, making relations with Algebraic Geometry more subtle. A famous scale from which to view such moduli spaces is the \textit{diameter scale}, for instance, in relation to the Strominger-Yau-Zaslow conjecture in Mirror Symmetry \cite{StromingerYauZaslow1996}. At such scale, however, natural moduli compactifications won't be algebraic at the boundary, but ``tropical'' (that is, with a topological boundary parametrizing certain \textit{real} spaces with affine structures). On the other hand, at \textit{volume scale} things seem to connect more to classical Algebraic Geometry over $\mathbb{C}$, pointing towards the existence of (differential) geometrically meaningful compactifications of moduli of polarized Calabi-Yau manifolds, that should be seen as analogous to the Baily-Borel compactification for polarized K3 surfaces (see Conjectural picture \ref{CP2}, and the subsequent discussion). This appears to be compatible with the recent works by, e.g.,  Sun-Zhang \cite{SunZhang2019}, Y. Li \cite{Li2023} and Odaka. At even higher scales in these collapsing situations one will see the appearance of possibly collapsing at infinity bubbling limits (e.g., certain 2D gravitational instantons \cite{SunZhang2024}, Tian-Yau metrics \cite{tian1990complete}, Collins-Li \cite{CollinsLi2024}, etc...), whose general relation with algebraic moduli compactifications also remains quite mysterious (see some related conjectures in \cite{Odaka2022b}). However,  I won't discuss these higher scales too much in this survey (but see remark \ref{HS}).\\

The structure of the survey is as follow. In the Section \ref{S1}, I describe a simple  motivational example, namely the moduli of flat metrics on the sphere. In such example one can basically fully exploit the multiscale aspects of the moduli spaces. In particular, this could motivate the conjectural pictures proposed later. I begin Section \ref{S2} with a discussion of non-collapsing bubbling in complex dimension two, and then, after recalling how to algebro-geometrically compute the metric tangent cones, I focus on AC bubbling in general dimensions, giving some examples and discussing the Conjectural picture \ref{CP1}. In the first part of Section \ref{S3}, I touch upon certain aspects of low scales (namely, the diameter and volume scales) compactifications of moduli of polarized KE varieties. In particular, I recall K-moduli spaces of varieties, and discuss a possible generalization in the full collapsing Calabi-Yau situation, aiming to provide some \textit{unifying} geometric viewpoint for K-moduli (see Conjectural picture \ref{CP3}). Finally, in the last part of Section \ref{S3}, I discuss related aspects, motivated by the Calabi-Yau picture, for collapsing of negative KE metrics (near isolated log canonical singularities), leading to the picture described in Conjecture \ref{CP3}.\\

I hope this survey can help young people in the field to have an overview of some of the current emerging new exciting directions, as well as of the general motivational picture behind them.\\

\textit{Acknowledgments} I would like to deeply thanks M. de Borbon, Y. Fang, Y. Odaka, S. Sun. for the many discussions and work we have done together and which are the basis for this survey. This survey is basically a written-up version of a talk firstly given during the conference \textit{Complex Hermitian Geometry}, 19-23 May 2025 Angers (France), celebrating Prof. Paul Gauduchon '80s birthday. \\

Finally, I like to thanks the Villum foundation for supporting with the grant Villum YIP+ 00053062.

\section{A simple example: moduli of flat metrics on the sphere} \label{S1}

A quite elementary, but nevertheless rich, example I'd like to illustrate is the one of moduli spaces of flat metrics with cone angles less than $2\pi$ (a convexity condition) on the two dimensional sphere, in particular in relation to multiscale aspects and their interaction with Algebraic Geometry. Most of the results in this section are taken from a series of works done in collaboration with de Borbon \cite{deBorbonSpotti2019, deBorbonSpotti2023, de2024some}. \\

The underlying geometric PDE in such examples is, of course, just:
$$K(g)=0,$$
with $K$ the Gauss curvature for $g$ metrics with cone angle singularities on $\mathbb{S}^2$. By cone angle singularities, we mean that such flat metrics $g$ are locally isometric to a 2D flat cone of angle $2\pi \beta < 2\pi$ near the tip:
$$ g=_{loc} dr^2 +\beta^2 r^2 d\theta^2.$$ 

Examples of such flat metric spaces are very concrete. For instance, the surfaces of convex 3D polyhedra (note that, from a metric view point, the edges do not give singularities, while the vertices are the cone points), or doublings of convex planar polygonal figures. The collection of solutions of such non-linear PDE form a non-trivial moduli space (up to scaling and diffeomorphisms), say $\mathcal{M}$.\\

In order to study such metrics from a more algebraic point of view, it is convenient to fix some parameters (markings), for instance, the values of cone angles $2 \pi \beta_i<2 \pi$ at the vertices $p_i$, and assume that the angles $\beta_i$  are rational (irrational ones can be similarly handled from a more general Symplectic viewpoint). After taking the natural conformal structure such flat metrics induce on the sphere, we can think at these shapes as the Riemann sphere $\mathbb{P}^1$ with $k$-marked points $\{p_1, \dots, p_k\}$ (determined up to projective transformations), all with attached cone points $2\pi \beta_i$. If we would fix the position of the cone points, the equation would then linearize to just a Laplacian on the conformal factor but, actually, the solutions are even completely explicit! Such flat metrics take the form
$$
g=c\left(\prod_{i=1}^{k-1} |z-p_i|^{2\beta_i-2}\right) dzd\bar z,
$$
in terms of standard complex variable $z$ on $\C$. Here one of the point is assumed, without lost of generality, to be equal to $\infty$. By the conical version of Gauss-Bonnet, we know that the sum of the cone angles must be equal to $\sum_i (1-\beta_i) =2$. Such unique up to scaling flat metrics \textit{geometrize} the algebraic log pairs 
$$(\mathbb{P}^1, \sum_i(1-\beta_i)p_i),$$ 
and it is then natural to expect they will be relevant in studying their geometry.

 \subsection{Moduli of flat metrics at the fixed diameter scale}

Fixing the values of the cone angles (or weights), up to scaling and holomorphic isometries such metric spaces form moduli spaces $\mathcal{M}(\underline{\beta})$ of complex dimension $k-3$. How can we describe the global structure of such moduli?  

First note that such moduli spaces are non-compact. For generic $\beta$s (more precisely, if $\sum (1-\beta_{i_j})\neq 1$), such metrics do \textit{not collapse} at fixed diameter scale. However, points can still \textit{collide} to each other, resulting in new flat metrics on the sphere with less, but sharper, cone angles. From a metric point of view it is then  very natural to \textit{compactify} the moduli spaces by adding such collisions as well. The topology considered here is essentially the  Gromov-Hausdorff (GH) topology, where we want also to recall the complex structure (this is possible since we have convergence up to diffeomorphism on the complement of arbitrarily small neighborhood of the collisions). Denote such compactification as $\overline{\mathcal{M}(\underline{\beta})}^{GH}$. 

It is not too hard to see that holomorphic isometries between two spaces linearize and then, thanks to the explicit description of the metrics in combination with some works by Deligne and Mostow \cite{DeligneMostow1986}, the compactified  moduli spaces $\overline{\mathcal{M}(\underline{\beta})}^{GH}$ can be described \textit{algebro-geometrically} as a concrete GIT quotient, with natural linearization induced by the fixed values of cone angles:

\begin{theorem}[\cite{deBorbonSpotti2019}, \textit{after} Deligne-Mostow \cite{DeligneMostow1986}]
There is a homeomorphism

$$\overline{\mathcal{M}(\underline{\beta})}^{GH} \cong (\mathbb{P}^1\times \dots \times \mathbb{P}^1)//_{\underline{\beta}} SL(2)=: \overline{M}^{GIT}_{\underline{\beta}},$$

with boundary strata corresponding to collisions of points.\\
\end{theorem}

We should remark here that this is one of the simplest example of $K$-moduli spaces of varieties (see Section \ref{S3} for a more general discussion on the topic): namely this moduli space is the coarse moduli spaces of the log  pairs $(\mathbb{P}^1, \sum (1-\beta_{i}) p_i) $, $1$-dimensional (K-stable) Calabi-Yau pairs, with its compactification as K-moduli.\\

If instead we assume the \textit{non-generic} condition $\sum (1-\beta_{i_j})=1$ (from the point of view of doublings, this corresponds to have some parallel edges), then such flat metrics \textit{collapse} to lower dimensional unit real intervals at fixed at diameter equal to one scale. Such collapsing to intervals happens precisely in correspondence of the isolated GIT polystable but not stable points. Thus, we can still give some \textit{topological }identification $$\overline{\mathcal{M}(\underline{\beta})}^{GH} \cong_{top} \overline{M}^{GIT}_{\underline{\beta}},$$but now with cuspidal points parametrizing real intervals. 

The simplest example of this phenomenon is the moduli space given by weights $(\frac{1}{2},\frac{1}{2},\frac{1}{2},\frac{1}{2})$. Such case is related to the moduli space of flat metrics on the torus by quotienting with the standard elliptic involution: the moduli spaces before compactification is isomorphic to $\mathbb{P}^1 \setminus \{0,1,\infty\}$, that is, as will be more clear in a moment, with the  $\infty$-pair of pants with three cusps. Forgetting the marking, we then recover the universal modular curve $\mathbb{H}/SL(2,\mathbb{Z})$,  that is, the moduli space of complex one dimensional tori. Examples of more complex Shimura-type moduli spaces also arise when considering more than $4$ cone points.

\subsection{Moduli of flat metrics: non-collapsing bubblings}

The above concludes the description of such moduli spaces at the fixed diameter scale. But what if we rescale the metric as points collide to each other? 

Let first consider \textit{rescalings} in the pointed GH topology $$\lim_{t\rightarrow \partial \overline{\mathcal{M} (\underline{\beta})}^{GH}} (\mathbb{P}^1, \lambda^2(t) g_t, p(t)),$$ as $\lambda(t)\rightarrow \infty$ and $t\rightarrow \partial \overline{\mathcal{M} (\underline{\beta})}^{GH}$ at collisions \textit{away} from the cuspidal regions (so that collapsing is not happening). The next scale we  see is the one of \textit{metric tangent cones at the singularities}. In this case just ordinary $2D$ flat cones of angles $2\pi\alpha$. Such angle can be computed from the collision of points via the formula $$1-\alpha=\sum_j (1-\beta_{i_j}).$$

If we rescale even further, we see the appearance of the first \textit{conical bubbles} (minimal bubble), that is, of conical metrics on $\C$ with conical geometry at infinity and residual conical singular points. In particular, the first cone at infinity is isometric to the cone at the singularities, say with cone angle $2\pi\alpha_0$. Such bubbles may have still conical singularities caused by residual collisions of points at such higher scale. Thus we can iterate the above scaling procedure to find new tangent cones at singularities (possibly changing the observer, i.e., the points we focus the rescaled GH limits) and thus corresponding conical bubbles. Since less and less points are colliding at each steps, the various cone angles $\alpha_i$ strictly increase until we reach in a \textit{finite} number of steps the deepest bubbles, namely rescaled limits for which all further rescaling of the degenerating sequence of flat metrics will be just tangent cones. \\

Such rescalings can be collected in the  \textit{metric bubbles' trees}: $$\tau_p:=\left\{\lim_{t\rightarrow \partial \overline{\mathcal{M} (\underline{\beta})}^{GH}} (\mathbb{P}^1, \lambda^2(t) g_t, p(t)), \, \mbox{with} \, p(t)\rightarrow p \in X_{\infty} \mbox{ for } \lambda(t)\rightarrow \infty\right\}/\mbox{isom}$$ for each collision cone point $p \in X_{\infty}\in\overline{\mathcal{M} (\underline{\beta})}^{GH}$. This object is indeed a connected tree made of tangent cones and conical bubbles. \\

The natural question which now arises is: can we compute such bubble trees $\tau_p$ \textit{algebraically}, that is, from the actual algebraic data of degenerating families? The answer is yes, and it can be easily be read from the equations of the colliding points. 

\begin{theorem}[\cite{de2024some}]\label{1DB}

Let $\pi:(\mathcal{X}=\mathbb{P}^1\times \Delta, \mathcal{D})\rightarrow \Delta$, with $\mathcal{D}=\sum_{i=1}^{k}(1-\beta_i)p_i(t)$ polynomials and $\sum_{i=1}^{k}(1-\beta_i)=2$. W.l.o.g., let $S=\{p_{i_j}\}$ be the subset of polynomials vanishing at zero $p_{i_j}(0)=0$ (collision) and equip the fibers with the unique volume $1$ conical flat metric. Then: 

\begin{itemize}
\item The metric bubble tree $\tau_0$ (of rescaled pGH limits) at zero can be computed algebro-geometrically from the polynomials by a sequence of nested relations on $S$: for $l\geq 0$,
$$ p_{i_j} \sim_l p_{i_k} \; \Leftrightarrow \;  ord_0(p_{i_j}-p_{i_k})\geq l, $$
where to each level $l$ equivalence class $S^l_m$ we associate a \textit{conical bubble} (unique up to scaling) whose conical points are determined by the value in zero of $l+1$ derivatives of polynomials in $S^l_m$.\\

\item A choice of holomorphic section $\sigma$ for $\pi$ with $\sigma(0)=0$ (the observer) naturally determines unique equivalence classes at various level, giving a \textit{path} from the root of the tree (minimal bubble) to a deepest bubble relative to $\sigma$ (not necessarily one of the leaves of the bubble tree). 

\end{itemize}
\end{theorem}

Let me point out that, even if the proof is just based on comparing the above algorithm with the explicit description of the metrics, the actual meaning of the above theorem is that, in principle, the bubbling is visible and computable only from the data of the algebraic family. This is precisely the type of algebraic computability of the trees we are looking for, although in an elementary situation.\\

Here a schematic picture of a bubbles tree: 

\begin{center}
   \includegraphics[scale=0.25]{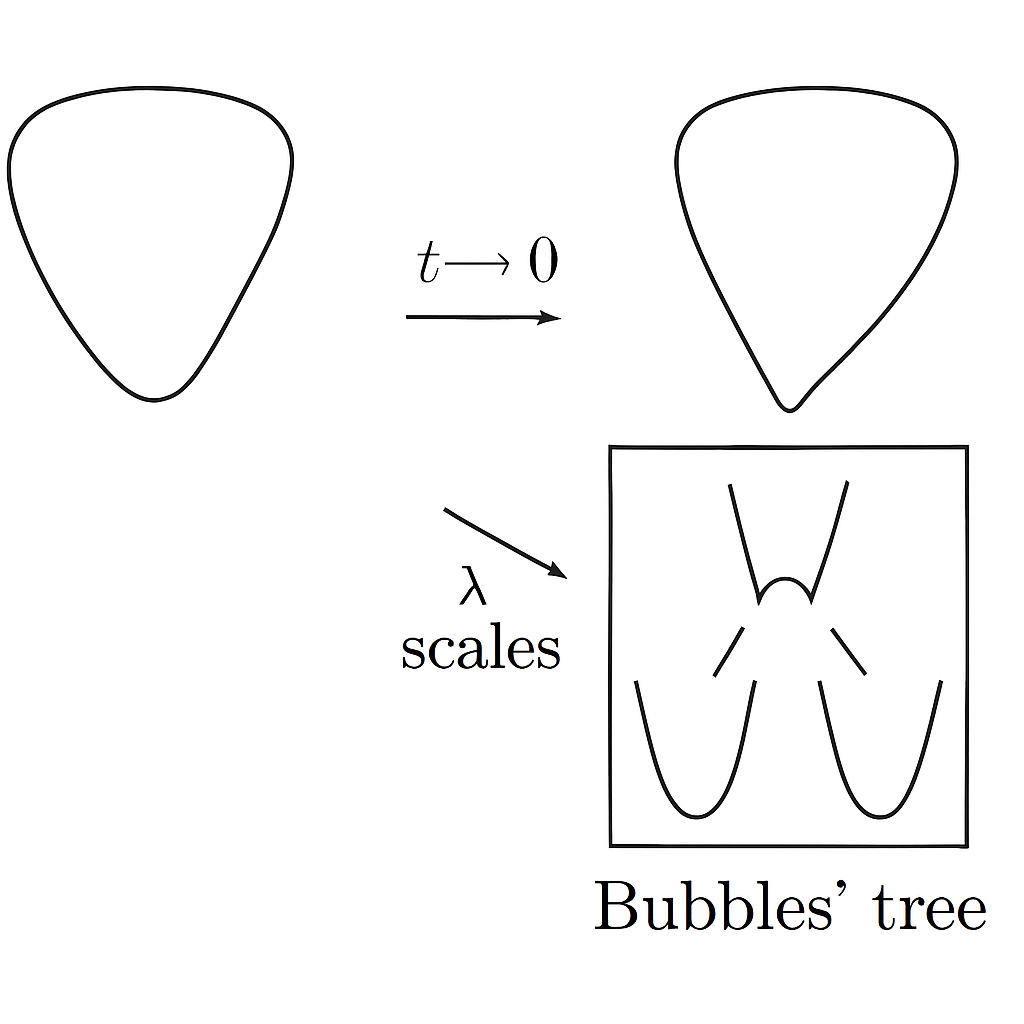}
\end{center}

So far we have investigated the multiscale geometry of the (non-collapsing) objects parametrized by the moduli spaces, but one can further ask about \textit{multiscale aspects of the coarse moduli space as well}. Indeed, such moduli $\mathcal{M}(\underline{\beta})$ spaces admit a natural metric, the so-called Weil-Petersson metric 
$$\omega_{WP}=i \bar\partial\partial \; log \, \left( \int \Omega_{\underline{\beta},t} \wedge  \overline{\Omega}_{\underline{\beta},t}\right), $$ where $\Omega_{\underline{\beta},t}$ is the locally defined meromorphic $1$-form with poles at the cone points. 

A remarkable fact (and one of the reasons for the interest of Deligne and Mostow for such moduli problem) is that such metric turns out to be complex hyperbolic, in particular, itself K\"ahler-Einstein with negative Einstein constant. This metric is incomplete. However, Thurston showed in \cite{Thurston1998} that near the non-collapsing region, its metric completion can be identified again with the GIT quotient $\overline{M}^{GIT}_{\underline{\beta}}$. 

From a complex differential geometric picture, one can now ask what the next scale/asymptotic of such WP metric near its completion is in terms of complex geometric data. The answer to that question can be given via a very standard construction of metrics on algebraic cones, namely the Calabi's Ansatz. At a generic point in $\partial \overline{M}^{GIT}_{\underline{\beta}}$, the completion gets a natural weighted log structure $(\overline{M}^{GIT}_{\underline{\beta}}, \partial \overline{M}^{GIT}_{\underline{\beta}})$, with $$\partial \overline{M}^{GIT}_{\underline{\beta}}=\sum_{i<j} (1-{\mu_{ij}}) D_{ij},$$ for coefficients $\mu_{ij}$ computed from the collision of just two points (generic boundary condition), via the formula $1-{\mu_{ij}}=2-(\beta_i+\beta_j)$. At each deeper non-collapsing point the pairs $\partial \overline{M}^{GIT}_{\underline{\beta}}$ has still Kawamata-log-terminal (klt) singularities, as a pair. The (log) algebraic structure of such singularities is given by product of pull-backs to $\mathbb{C}^{j_i+1}$ of some special (weighted) hyperplane arrangements in $\mathbb{P}^{j_i}$ times $\mathbb{C}^{k-3-\sum_i(j_i+1)}$. By comparing with Thurston's results (see also \cite{CouwenbergHeckmanLooijenga2005}), it is then possible to check that the metric structure is indeed compatible with such description. First, at generic points in the non-collapsing boundary, the metric is asymptotic to a cone of angle $2\pi \mu_{ij}$  times flat $\mathbb{C}^{k-4}$. At deeper points (i.e., not generic in codimension $1$) the WP geometry is then described by:

\begin{theorem}[\cite{deBorbonSpotti2019}, \textit{after} Thurston \cite{Thurston1998}] 
 
 Near a deeper non-collapsing boundary point the WP metric has metric tangent cone given by a \emph{Calabi's ansatz} cone construction from products of constant positive holomorphic sectional curvature  on braid hyperplane arrangements $((\mathbb{P}^{j_i},\underline{\mu}H), g_{KE})$ times flat factors, i.e., $$g_{WP} \cong dr^2+r^2\overline{g},$$ 
 for small $r$, where $\overline{g}$ is the standard Sasaki lift of  $g_{KE}$ to a metric on the sphere $\mathbb{S}^{2j-1}$, conically singular along circle bundles in the Hopf's fibration.

\end{theorem}

Before discussing the case of collapsing, let's remark some further relations between the asymptotic of the WP metric and the actual conical bubbling. First of all, note that the complex link factors $\mathbb{P}^{j}$s above, with $j=\sharp S - 2$, can be identified with the moduli spaces of minimal bubbles after factoring out by translations and rescalings. More generally, it is possible to observe that the Deligne-Mumford compactification (a log resolution of the moduli spaces) $$\overline{\mathcal{M}}_{0,k}^{DM}\rightarrow \overline{M}^{GIT}_{\underline{\beta}}  \cong \overline{\mathcal{M}(\underline{\beta})}^{GH}$$ has a further differential geometric meaning, different from the usual one in term of cuspidal (at the marked points) hyperbolic metrics: by identify the \textit{main component} of a DM-stable degeneration with the diameter scale limits, and the further components  as bubbles (intersecting if the tangent cone at a singularity matches the tangent cone at infinity of a bubble), one can give to the DM compactification  the differential geometric meaning of moduli spaces of bubbles' trees. Hence, the Deligne-Mumford compactification can be viewed also as a ``resolution of all scales'' of the degenerating flat metrics. For more discussion on this, see \cite{de2024some}.

\subsection{Moduli of flat metrics: collapsing}

Next we consider scalings  $(\mathbb{P}^1, \lambda^2(t) g_t)$ as $\lambda(t)\rightarrow \infty$ and $t\rightarrow \partial \overline{\mathcal{M} (\underline{\beta})}^{GH}$ at the cuspidal regions, i.e., when collapse happens. Here we see the following multiscale limits:

\begin{itemize}

\item \textit{diameter scale}: convergence to an \textit{interval} $I$ (``large complex structure equal to small complex structure limit'');

\item \textit{volume scale}: generically convergence to the real line  $\mathbb{R}$ (``metric tangent space to I''), or half lines.

\item \textit{first non-collapsing scale}: convergence to a bubbling  \textit{cylindrical} cones metrics on $\mathbb{C}$ (``Tian-Yau type metrics''), or a cylinder $S^1\times \mathbb{R}$.

\item \textit{higher scale}s: bubblings in conical metrics on $\mathbb{C}$, as in non-collapsing case.

\end{itemize}

In the above parenthesis I used the terminology which is nowadays common for describing general Calabi-Yau degenerations. Note that at lower scales we see indeed collapsing to lower dimensional metric spaces. \\

Also in this collapsing situation one can analyze the behavior of the WP metric on the moduli space, which is now a complete metric near the cusps. 

\begin{theorem}[\cite{deBorbonSpotti2019}, \textit{after} Thurston \cite{Thurston1998}]

Near collapsing cusps the WP metric is complete and isometric to a \emph{Calabi's ansatz} cone construction from some flat hyperplane arrangements $((\mathbb{P}^{k-4}, \underline{\mu}H), g_{flat})$ (or product of two):

$$g_{WP}=_{loc}g_{cusp} = \frac{1}{2}dt^2 + e^{-t}(\pi^*g_{flat} + 2e^{-t}d\theta^2), \mbox{ and } \omega_{cusp} =i\bar\partial\partial t.$$
\end{theorem}

In such situation the cuspidal WP metric is, of course, not (close to) a Riemannian cone and indeed itself manifests interesting multiscale collapsing: if we do not rescale, the metric fully collapses to the real line  $$g_{cusp}\rightarrow_{pGH} \mathbb{R},$$ as $t_i\rightarrow \infty$. However, there is another interesting scale where we see the emergence of the flat  Calabi-Yau geometry via the collapsing of the \textit{horospheres}  $t=c$ (that is, the level set of the natural K\"ahler potential of the WP metric), i.e., $$e^{t_i} g_{cusp} \rightarrow_{pGH} \mathbb{R}\times ((\mathbb{P}^{k-4}, \underline{\mu}H), g_{flat}).$$ 

Note also that $$K_{\overline{M}^{GIT}_{\underline{\beta}}}+ \partial \overline{M}^{GIT}_{\underline{\beta}}=K_{\overline{M}^{GIT}_{\underline{\beta}}}+\sum_{ij} (1-{\mu_{ij}}) D_{ij}>0$$ is ample, and $(\overline{M}^{GIT}_{\underline{\beta}},\partial \overline{M}^{GIT}_{\underline{\beta}})$ a log canonical pair. Thus $\omega_{WP}$ is an example of negative KE metric with cone singularities along divisors, showing such interesting multiscale behavior near its cusps, precisely the log canonical (lc) singularities of the pair.

Finally, similarly to the non-collapsing case, one can interpret such splitting Calabi-Yau flat limits of the horospheres as moduli spaces of minimal non-collapsing bubbles. Indeed, the actual GIT quotient $\overline{M}^{GIT}_{\underline{\beta}}$ is singular at the cusps with singularities given by cone over products of two projective spaces $\mathbb{P}^{k-4-j}\times \mathbb{P}^{j}$, where each of the two factors parametrizes cylindrical collapsing bubbles for the exactly two regions of collision.

From an algebraic perspective, in order to see such bubbling it is quite significant to perform some elementary birational modifications of an original product family $\mathcal{X}=\mathbb{P}^1\times \Delta$ analogous to Theorem \ref{1DB}, now for polynomials sections clustering to zero with $\sum (1-\beta_{i_j})=1$ and the remaining ones also going to infinity (and co-angles summing to one): blow-up  the total space of the family at the zero and infinity clustering points of the central fibers. Then we have a new family $\hat{\mathcal{X}}$ with three components $E_0, \overline{\mathbb{P}^1}, E_1$ in the central fiber where the exceptional $E_i$s can be identified with the cylindrical bubbles (with cone points determined by strict transforms of the polynomial sections), while the central one $\overline{\mathbb{P}^1}$ as a connecting cylinder, which you can get rid off by contracting it and creating a further new family $\tilde{\mathcal{X}}$ with $A_1$-singularity located at the central fiber. Note that the dual complex of the central fiber (that is, replace the two components with points, and connect them with a line, since they intersect) of this minima model is an interval (as the diameter scale limit), with the volume scale geometry focuses near the intersection (log canonical center) of the two higher scale bubbling components. Of course, this is very similar to the discussion about non-collapsing bubbling in relation to the Deligne-Mumford compactification we discussed before.

\subsection{Further considerations}

Before concluding this first section, let me mention few more facts, which are also somewhat related to multiscale geometry.\\

1.\textit{ Generalized BMY inequality}. These spaces $S=(X, \underline{\mu}H)$ of constant holomorphic sectional curvature emerging at various scales (the actual moduli space or the Fano/Calabi-Yau hyperplanes' arrangements appearing at its singularities) are rather special, for instance they satisfy equality in the  \textit{Bogomolov-Miyaoka-Yau} (BMY) inequality, provided  codimension two corrections related to metric tangent cones are added to the second Chern class $c_2$. For instance, in dimension two, one defines a new $c_2$ via the following formula:
$$c_2^\nu(S)= e(X)+\sum_i (\mu_i-1) e(H_i \setminus Sing(H) )+\sum_{p\in Sing(H)} (\frac{1}{4}\hat{vol}(\nu_{min}(p))-1)$$
 where $\hat{vol}(\nu_{min}(p)$ is the Li's normalized volume (see later Section 3.3), a purely algebraic invariant of the singularities of the pair, which coincides with the metric density of the tangent cones $\hat{vol}(\nu_{min}(p))= 4 \Theta_p$, and $e$ the usual topological Euler characteristic. It is then natural to expect such correction to appear in BMY for pairs in general, as it has been discussed in some two dimensional cases in \cite{deBorbonSpotti2023} (see also Li \cite{Li2021} and Langer \cite{Langer2003}).\\

 2. \textit{Modular interpolation}. Some interesting observations arise when we note that these moduli spaces of flat metrics $\overline{\mathcal{M}(\underline{\beta})}^{GH}$ sit between  other moduli spaces of conical metric with positive and negative Gauss curvatures
 $$K(g)=\lambda,$$
 with $\lambda \in \mathbb{R}.$ Indeed, if the condition on the sum of the co-angles is changed to $$\sum_i(1-\beta_i)=2 \pm \epsilon,$$ the pair $(\mathbb{P}^1, \sum_i(1-\beta_i)p_i)$ becomes Fano (so $\lambda>0$) or of general type (so $\lambda<0$). 
 
 For negative epsilon (Fano situation), the GH compactification will also agree with the GIT quotient, now with polystable points parametrising (at volume$=$ diameter scale) ``rugby balls'', that is positively curved conical metrics on the sphere with two equal angles $\alpha_0=\alpha_{\infty}$ at zero and at infinity. These are the simplest examples of K-moduli spaces of KE Fano (pairs). 
 
 For positive epsilon (general type situation), the natural moduli spaces compactifications are the Hassett moduli spaces \cite{Hassett2003}, which are the easiest example of general KSBA moduli, still K-moduli for positive canonical class. In such cases, however, the moduli space will go through birational changes starting at polystable points, essentially analogous to the ones mentioned before to describe cylindrical collapsing. From a metric perspective, such components  carry (at volume scale) complete conical hyperbolic metrics with cusps at the intersection of the divisors. 
 
 If we look at the metric geometry at volume scale as we pass through the Calabi-Yau \textit{threshold } $\sum_i(1-\beta_i)=2$, we see from the Fano direction the rugby balls stretching and converging to a real line, while from the negative curvature side we see along the degeneration process  that the volume starts concentrating more and more near the cuspidal formation region, leading in the limit still to a convergence to a line. The cylindrical bubbling will be visible at higher scales (as well as further conical bubblings). A similar picture holds at the diameter scale, with convergence to an interval.

A further final remark. Such continuity aspects have been exploited for computing volumes of such moduli spaces by Tambasco \cite{Tambasco2022}, who recovered McMullen's computation  \cite{McMullen2017} of the $\omega_{WP}$ volume in the Calabi-Yau case as a limiting of CM bundle degrees from the Fano direction.

 \begin{remark} There is also a more differential-geometric/analytic perspective on moduli spaces of conical metrics on surfaces (including bubbling) recently developed by Mazzeo and Zhu \cite{MazzeoZhu2020}.
 \end{remark}

\section{Non-collapsing: asymptotically conical bubbling} \label{S2}

We now move to the general case of Kähler-Einstein metrics on some $n$-dimensional complex manifolds $X_t$. Thus, the geometric PDE we consider is, for Einstein constant $\lambda \in \mathbb{R}$,
$$Ric(\omega_t)=\lambda \omega_t,$$
where $\omega_t$ is the K\"ahler form and $$Ric(\omega_t)=i\bar\partial \partial \log (\omega_t^n),$$ the Ricci form. For us $t$ denotes a complex structure parameter, and we assume that the KE metrics are \textit{polarized}, i.e., $[\omega_t]\in c_1(L_t)=\alpha$ fixed. For non-zero $\lambda$ constant, the polarization is just given by plus or minus the canonical bundle $L_t=K_{X_t}^{\pm 1}$. 

Such KE metrics on $(X_t, L_t)$ are unique (up to automorphism), hence they geometrize the polarized manifolds, in full analogy with the complex one dimensional example described previously. However, in this higher dimensional situation, when the complex structure and polarization are fixed, the KE equation does not linearize even in the $\lambda=0$ case, but instead it reduces to the scalar Monge-Amp\`ere equation (first solved in the famous $\lambda=0$ Calabi-Yau case by Yau in \cite{Yau1978}, but with solutions, in general, absolutely not explicit).\\

In this section, we focus on the study of such KE degenerations in the \textit{non-collapsing} case. This means that for a given family $\omega_t$ of polarized KE metrics (with same Einstein constant), we assume that we have the following uniform control on small balls:
$$Vol(B_p(r))\geq \varepsilon r^{2n},$$
for some $\varepsilon$ independent of $t$. In the case of positve Einstein constant (i.e., Fano case), this condition is always automatically satisfied (depending only on the degree of the Fano), thanks to Myers and Bishop-Gromov theorems. For $\lambda \leq 0$ (i.e., Calabi-Yau or ample canonical bundle cases), this is an extra condition to impose. Thanks to Donaldson-Sun theory \cite{donaldson_sun_gh_limits_2014},  Gromov-Hausdorff limits of such non-collapsing polarized KE manifolds are going to be some singular algebraic varieties $X_0$ with at most klt singularities, and which are naturally homeomorphic to the metric completion of the smooth limit KE metric on $X_0\setminus Sing(X_0)$. \\

Understanding such compact limits $X_0$ has been relevant in the study of the KE/K-moduli compactifications of moduli spaces (as I will recall in more details at the beginning of the next Section \ref{S3}). However, the problem we want to analyze here is the following: under such non-collapsing conditions, what can the \textit{rescaled bubbling limits} (in the pointed Gromov-Hausdorff topology) $$\lim_{t\rightarrow 0,  \lambda(t)\rightarrow \infty} \left(X_t, \lambda^2(t)\omega_t, p(t)\right)=B_\infty$$ of the degenerating KE metrics be, and how are they possibly related to Algebraic Geometry?

\subsection{The complex two dimensional case}

Let's begin with the complex two dimensional case. From a differential geometric/analytic point of view the bubbling theory is very well-known and it goes back to the seminal works by Anderson, Bando-Kazue-Nakajima for non-collapsing limits of general Einstein $4$-manifolds in late $'80$s.

\begin{theorem}[e.g., \cite{Anderson1989, BandoKasueNakajima1989}] Given a non-collapsing sequence of KE surfaces $(X_i,\omega_i$), its rescalings lead to finite bubbles trees of quotients of orbifolds Kronheimer's ALE spaces.

\end{theorem}

Recall that Kronheimer's ALE spaces \cite{Kronheimer1989} are the non-compact hyperk\"ahler manifolds which away from a compact set look like $\C^2/\Gamma$ with $\Gamma$ a finite subgroup of $SU(2)$, and with metric asymptotic to the flat one. Such metrics are explicitly described in terms of hyperk\"ahler quotients from  a finite dimensional flat space and connect to the deformation theory of ADE singularities, e.g., $A_k$-singularities $$z_1^2+z_2^2=z_3^{k+1}.$$ The simplest example for such manifolds is the so-called Eguchi-Hanson space $$X\cong T^*S^2 \cong \{z_1^2+z_2^2+z_3^2=1\}\subset \mathbb{C}^3$$ with its very explicit ALE metric $$\omega=(i\partial\bar\partial \sqrt{\vert z \vert^2+1})_{\vert X}.$$

\begin{remark} In general, one needs to consider ALE spaces with residual orbifold singularities (i.e., they have local singularities of type $\mathbb{C}^2/\Gamma_p$ (with $\Gamma_p \subset SU(2)$ of smaller order compared to the group $\Gamma$ at infinity), which are still capture by Kronheimer's construction. Moreover, one needs to consider \textit{quotients} of Kronheimer's ALE (e.g., compare Suvaina \cite{Suvaina2012}), with torsion canonical bundle. Indeed,  such \textit{torsion bubbles} can indeed appear as rescaled limits of Einstein metrics, as it has been shown in \cite{odaka2016compact} for degenerations of certain positive KE (del Pezzo) surfaces.
\end{remark}

The next natural question to consider is how to compute such bubble trees \textit{{purely algebraically}} given a degenerating polarized family of KE  varieties $X_t$ (equivalently, thanks to the YTD conjecture, K-polystable ones): $$\mathcal{X}\rightarrow \Delta,$$ with $X_0$ with, say, ADE singularities.

A local(!) answer to this question can be given by the next result, that generalizes Theorem \ref{1DB} to the complex two dimensional $A_k$-case. Consider the following explicit hyperkähler ALE metrics given by the \textit{Gibbons-Hawkin}g's ansatz \cite{GibbonsHawking1978}:
$$g(t)=V(t)(dx_1^2+dx_2^2+dx_3^2) + V^{-1}(t) \eta^2$$
with  $V$ the harmonic function on flat $\mathbb{R}^2$ $$V(t)=\frac{1}{2}\sum_{i=1}^{k+1} \frac{1}{\vert x-x^i(t)\vert},$$ with polynomials $x^i(t)=(x_1^i(t),x_2^i(t), 0)$ vanishing at $t=0$, and $\eta$ the connection one form  for an $S^1$-bundle (of Chern number one) over the punctured plane, such that $$d\eta=-*d V(t).$$

It turns out that the above metrics $g(t)$ extend to define for each $t$ a $i\partial \bar \partial$-exact Kähler metric on the affine complex manifold:
$$z_1^2+z_2^2=\prod_{i=1}^{k+1}(z_3-p_i(t)),$$
where $p_i(t)=x_1^i(t)+\sqrt{-1}x_2^i(t) \in \C$. The case $k=1$ gives the Eguchi-Hanson space described before. Note that that affine manifolds gives a curve in the (Galois cover) of the versal deformation space of an $A_k$-singularity, thus such ansatz provides a fully general description of what happens locally near such singularity. 

As in the two dimensional flat case case, we want to compute the bubbles' tree for these metrics $g(t)$ as $t \rightarrow 0$. The following theorem, based on explicit computations using the Gibbons-Hawking's ansatz description of the metric, provides a complete answer to this problem.

\begin{theorem}[\cite{deBorbonSpotti2023}] Analogously to Theorem \ref{1DB}, for the family of metrics $g(t)$ defined above on the varieties $$z_1^2+z_2^2=\prod_{i=1}^{k+1}(z_3-p_i(t)),$$ one can compute algebro-geometrically the bubble tree $\tau_0$ made of orbifolds Kronheimers' ALE spaces from the set of polynomial  $S= \{p_1(t), . . . , p_{k+1}(t)\}$ by grouping them according to their relative order of vanishing at the origin. \\

Moreover, by picking a section for the family we can determine a path from the tangent cone to a (relative to the section) deepest bubble. Such path of bubbles can be computed via a series of \emph{iterated rescalings} of the equation (after centering the rescaling points to zero) of the form:
$$z^j_1=t^{\lambda_j (k_j+1)}z^{j+1}_1, \; z^j_2=t^{\lambda_j (k_j+1)}z^{j+1}_2, \; z^j_3=t^{2 \lambda_j } z^{j+1}_3$$
as $t$ goes to zero for some $\lambda_j$ and $k_1=k > k_2...$. Such $k_j$ are the weights (corresponding to Li's minimazing valuation, see later Section 3.3) for a residual $A_{k_j}$ singularity that a (bubble) limit can have at the strict transform of the given section.

\end{theorem}

Let us describe how such rescaling works with an example. Consider a smoothing  of the $A_2$-singularity (equipped with the metric $g(t)$ given by the Gibbons-Hawking's ansatz) of the form $$x^2+y^2=z(z-t)(z-t^2),$$ 
and consider the section (observer) to be $\sigma(t)=(0,0,0)$. This is the simplest example of multiscale bubblings. The rescaling algorithm in such case goes as follow:\\

\begin{itemize}

\item  \textbf{First scale}: Using the weights $(3,3,2)$ corresponding to the $A_2$ singularity, consider $x=t^{3\lambda}x'$, $y=t^{3\lambda}y'$, $z=t^{2\lambda}z'$, and re-write  $$x'^2+y'^2=z'(z'-t^{1-2\lambda})(z'-t^{2-2\lambda}).$$ For $\lambda_1=\frac{1}{2}$ the family as a holomorphic limit, and we get $$x'^2+y'^2=z'(z'-1)(z'-t)$$ with central fiber \textit{jumping} to $x'^2+y'^2=z'^2(z'-1).$ By comparing the ansatz metric with such scaling operation,  \textit{this computes the \textbf{minimal} bubble: an ALE orbifold on $x'^2+y'^2=z'^2(z'-1)$, with tangent cone at infinity given by $\mathbb{C}^2/\mathbb{Z}_3$ but with a residual $A_1$-singularity}.\\

\item \textbf{Second scale:} Use now $A_1$ weights: $x'=t^{2\lambda} x''$, $y'=t^{2\lambda} y''$, $z'=t^{2\lambda} z''$. As before, for $\lambda_2=\frac{1}{2}$: $$x''^2+y''^2=z''(tz''-1)(z''-1),$$ whose central fiber is now a smoothing of $A_1$. \textit{This computes the \textbf{deepest} bubble, namely the Eguchi-Hanson space}.\\

\end{itemize}

If we would have instead chosen a deformation of type  $$x^2+y^2=z(z-t)(z-at),$$  we would have ended up in a bubble isometric to the metrics $g(t)$ themselves, since in such family the metrics indeed differ only up to scaling. Families of type $$x^2+y^2=z(z-t)(z-at+t^2)$$ would instead give a single (minimal equal to deepest) bubble on the ``tangent'' smooth ALE space $x^2+y^2=z(z-1)(z-a)$. The crucial aspect of the first above chosen family was that the path in $t$ such deformation defines in the versal deformation space is tangent to the walls parametrizing deformations with a residual $A_1$-singularity.

\begin{remark}
The above algebraic picture is expected to describe bubblings of general KE metrics also in the compact case, e.g., in combination with gluings  (see for instance \cite{BiquardRollin2014, spotti2014deformations} for a no-multiscale situation, or the results of Ozuch for general Einstein multiscale bubblings \cite{ozuch2022noncollapsed1, ozuch2022noncollapsed2}). For instance, it is natural to expect that KE metrics on K3 surfaces are close to such Gibbons-Hawking  model description near a non-collapsing degeneration. However, as we will see, we can also try to argue more a-priori, making the bubble trees to emerge from degenerations of compact KE manifolds  by linking directely abstract metric rescalings with the algebraic geometry of singularities' formation.
\end{remark}

\subsection{Bubbling in higher dimension}

For general dimension, the problem we want to discuss is then the following:

\begin{problem}Given a polarized family of KE (eq. K-polystable) varieties $\mathcal{X}\rightarrow \Delta$ with klt singularities, compute algebraically its bubbles trees.
\end{problem}

The first result we need is the following recent results of Sun, based on applications of H\"ormander's techniques refining previous works by Donaldson and Sun \cite{donaldson_sun_gh_limits_II_2017}, which provides an essential first link between the differential geometric rescaled limits of non-collapsing KE metrics and Algebraic Geometry.

\begin{theorem}[\cite{SSBubbling}] Let $(X_i, \omega_i)$  be a non-collapsing sequence  of polarized KE manifolds. Then:
\begin{enumerate}

\item It exists (up to subsequences) an \emph{affine  algebraic Calabi-Yau asymptotically conical (AC) minimal bubble} $(B, \omega_B)$ whose tangent cone at infinity is equal to the unique tangent  cone $C(Y)_p$ at a singularity of the Gromov-Hausdorff limit $X_0$. 

\item Taking the tangent cones to  singularities of the (minimal) bubble, one can iterate point ($1$) above to obtain deeper AC bubbles until, in a \emph{finite} number of steps, we reach deepest bubbles.
\end{enumerate}
\end{theorem}

 In such general higher dimensional situation the Calabi-Yau cone $C(Y)_p$, a singular Ricci-flat space describing the asymptotic geometry of the KE metric on $X_0$ near a singularity $p\in Sing(X_0)$, does not need to be locally biholomorphic to a neighborhood of the singularity in $X_0$ (``jumping phenomena'') nor to be with isolated singularities if $p\in Sing(X_0)$ is such. These features already happen in the case of three dimensional $A_k$ singularities for $k\geq 3$ $$z_1^2+z_2^2+z_3^2+z_4^{k+1}=0$$ where the metric tangent cone is biholomorphic to flat $\mathbb{C} \times \C^2/\mathbb{Z}_2$. We remark that, interestingly, such jumping phenomena do not occur in the more restricted hyperk\"ahler situation, thanks to very recent works by Namikawa and Odaka \cite{namikawa2025canonical}. Thus, such general hyperk\"ahler case should resemble a bit more the complex two dimensional situation previously described. The bubble $(B, \omega_B)$ has klt singularities and a maximal volume growth Ricci-flat metric. Being itself a limit of smooth Einstein metrics, it has unique Calabi-Yau tangent cones at its singularities \cite{donaldson_sun_gh_limits_II_2017}. Note that, contrarily to the complex two dimensional case, the bubble may not be diffeomorphic  away from a compact set to its tangent cone at infinity  $C(Y)_\infty= lim^{GH}_{\lambda \rightarrow 0} (B, \lambda \omega_B)$. \\
 
 A crucial point is that in the above ``a-priori'' results of Sun, both the tangent cones and the bubbles, despite being algebraic objects,  are not at all explicit (they are related to abstract algebraic limits of the given varieties induced by the unknown KE metrics). Thus: how can we identify/compute them concretely?

\subsection{Algebro-geometric computability of tangent cones}
Remarkably,  tangent cones are indeed already known to be computable algebraically. In the complex two dimensional  $A_k$ example above, the weights
$(k+1,k+1, 2)$  correspond to the natural scaling action of the cone metric on $\mathbb{C}^2/\Gamma_p$. To see that, take the orbifold chart $$z_1+iz_2=u^{k+1}, z_1-iz_2=v^{k+1}, z_3=uv. $$ Generalizing results by Martelli, Sparks and Yau \cite{martelli2006geometric} emerged in the context of the AdS-CFT correspondence Li has introduced in \cite{Li2018NormalizedVolumes} a function $\hat{vol}_p$ defined on the set (aka the non-Archimedean link of a singularity) of \textit{valuations} centered at $p \in Sing(X)$ that computes such scalings:
$$\hat{vol}_p(\nu)= A^n(\nu) vol_p (\nu),$$
where $vol (\nu)$ is a (local) Riemann-Roch quantity 

$$vol_p (\nu)=limsup_{k\rightarrow \infty}  \frac{\, dim (\mathcal{O}_{X,p}/\{f \, \vert \, \nu(f) \geq k\})}{k^n/n!},$$
with $\mathcal{O}_{X,p}$ space of germs of holomorphic function at $p\in X$, and $A(\nu)$ the log-discrepancy, e.g., for a divisorial valuation $ord_D(f)=ord_D (\pi^*f)$, $f \in \mathcal{O}_{X,p}$, for $E\subset \hat{X} \xrightarrow{\pi}  X$ on some birational model over $p$,
$$A(ord_E)=1+ord_E(K_{\hat{X}}-\pi^*{K_X})>0,$$
definition which extend by linearity for the ``dense'' set of quasi-monomial valuations $A(\nu)=\sum w_j A(ord_{E_j})$, with $w_j \in \R^+$ and $E_j$ normal crossing at $q\in \hat{X}$. 

For example, if we take the ordinary double point $$z_1^2+z_2^2+z_3^2=0,$$ and we consider $E$ the exceptional divisor in the minimal resolution, we  see that $A(ord_E)=1$ (since it is a crepant resolution) and $$vol_p (ord_E)= 2(=c_1(\mathcal{O}_{\mathbb{P}^1}(2))),$$ hence $\hat{vol}(ord_E)= 2$. More concretely, this valuation gives order two to the coordinate functions $z_i$, that is, it corresponds to the weight $\nu= (2,2,2)$, which is precisely the scaling of the flat metric on $\mathbb{C}^2/\mathbb{Z}_2$.

Using purely algebro-geometric considerations \cite{Blum2018Existence}, it has been proven that for any germ of klt singularity $\mathcal{V}_p$ the function $\hat{vol}_p$ has a unique (up to scaling) minimizer, hence it defines a new invariant for such singularities:
$$\hat{vol}(\mathcal{V}_p):=\inf_\nu \hat{vol}_p= \hat{vol}_p(\nu_{min}).
$$
In the above $A_1$ example, the valuation $\nu= (2,2,2)$ is indeed the minimizing one.

Coming back to the metric tangent cone, in \cite{LiWangXu2021TangentCones} it was proved that for non-collapsing limits of polarized KE metrics such minimizing valuation $\nu_{min}$ centered at the singularities indeed computes algebraically the metric tangent cones via  a $2$-step procedure, confirming a picture of conjectural picture by Donaldson-Sun \cite{donaldson_sun_gh_limits_II_2017}:
$$\mathcal{V}_0\subset X_0  \leadsto W=Spec(Gr(\nu_{min})))\leadsto_{test} C(Y),$$
Here $X_0$ is the non-collapsed GH limit, $W$ the semistable tangent cone, built from local holomorphic functions graded using the minimizing valuation, which degenerate along a one parameter subgroup to $C(Y)$ the $K$-polystable Calabi-Yau tangent cone with its Ricci flat cone metric. In such case, $$\hat{vol}(\mathcal{V}_p)=n^n \Theta_p,$$ where $\Theta_p$ denotes the \textit{metric volume density}
$$\Theta_p=\lim_{r\rightarrow 0}\frac{Vol(B_p(r))}{\omega_{2n} r^{2n}},$$
of the singular KE metric on $X_0$, which is also equal, up to some normalizing functor, to the volume of the link of the Calabi-Yau metric tangent cone.\\

Roughly speaking, the higher the normalized volume the less singular the klt singularity is. For instance, in dimension two $$\hat{vol}(\mathbb{C}^2/\Gamma)= 4/\vert \Gamma \vert.$$ In this direction, let me recall the following conjecture, which states that there should be a gap between the normalized volume of an ordinary double point (ODP) and the smooth case, where  $\hat{vol}(\mathbb{C}^n)=n^n$.
\begin{conjecture}[ODP Conjecture \cite{spotti2017explicit, liu2019k}] For a non-smooth $n$-dimensional klt singularity 
$$\hat{vol}(\mathcal{V}_p)\leq 2(n-1)^n  ,$$
with equality holding precisely for the ordinary double point singularity $\sum_i {z_i^2}=0$.
\end{conjecture}

Such conjecture has been shown to hold in dimension $3$ by Liu, Xu \cite{liu2019k} using Mori and Reid's classification of canonical singularities, and for general K-semistable cones \cite{LiMiao2025VolumeKSemistableFano}. It trivially holds also in dimension $2$, since  klt singularities are quotients and $\hat{vol}(\mathbb{C}^2/\Gamma))=4/\vert \Gamma \vert$. From a differential geometric point of view, it says that the Stenzel's Calabi-Yau cone is precisely the one with highest density among non-smooth ones, and hence it can be viewed as the KE analogous of the famous Willmore Conjecture, in its rephrasing in terms of minimal 3D cones. This type of concrete bounds on the singularities is relevant in the study of explicit moduli compactifications, compare \cite{odaka2016compact, spotti2017explicit, liu2019k}. Indeed, if true it would imply that the GH compactification of the moduli spaces of KE metrics on cubic hypersurfaces agrees with the natural associated GIT quotient (coming from the defining embedding) in \textit{every} dimension \cite{spotti2017explicit}.

More recently, Xu and Zhuang have proven \cite{xu2024boundedness} that there are \textit{no accumulations} away from zero of values of $\hat{vol}(\mathcal{V}_p)$ in each dimension. This can be viewed as an algebro-geometric counterpart to Sun's finiteness of bubbling: indeed, by Bishop-Gromov monotonicity, each deeper bubble would have tangent cone at infinity of lesser density (equal to normalized volume) compared  to its singularities. Since such density values do not accumulate, there will be only \textit{finitely} many different cones, hence interesting non-collapsing scales, describing the degenerations. \\

With these results in place, it is now time to formulate a general picture regarding the algebraic geometry of bubblings.

\subsection{Algebro-geometric computability of bubble trees}

In the complex dimensional two case, we have shown that bubbling is connected to the deformation theory of the quotient singularities  $\mathbb{C}^2/\Gamma$. However, in higher dimensions, the jumping phenomena of the tangent cones create some difficulties. A first issue is that even isolated singularities may have jumping tangent cones with infinite dimensional versal deformation space, as already the higher dimensional $A_k$ singularities show. For instance, the versal space of deformation space for the singularity $$\mathbb{C}^2/\mathbb{Z}_2\times \mathbb{C}=\{z_1^2+z_2^2+z_3^2=0\}\subset \mathbb{C}^4,$$ that is, the polystable tangent cone for $A_k$ singularity for $k\geq 3$, is \textit{infinite} dimensional, since we can add any polynomial in the $z_4$ variable: 
$$z_1^2+z_2^2+z_3^2=p(z_4)(=t_0+t_1z_4+t_2z_4^2+ \dots.)$$
However, if we consider \textit{negative weight deformations} relative to the weight action coming from the natural Calabi-Yau cone, or minimizing valuation, which in this case is $(2,2,2,1)$ (given by the product of the two flat metrics on $\mathbb{C}^2/\mathbb{Z}_2$ and $\C$), we see that
$$\lambda^4(z_1^2+z_2^2+z_3^2)=t_0+\lambda t_1 z_4+\lambda^2 t_2z_4^2+ \dots,$$
or, equivalently,
$$z_1^2+z_2^2+z_3^2=\lambda^{-4}t_0+\lambda^{-3} t_1 z_4+\lambda^{-2} t_2z_4^2+ \dots.$$
Thus the set of parameters $t_i \in Def(\mathbb{C}^2/\mathbb{Z}_2\times \mathbb{C}$) giving negative weights is only up to $t_3$, i.e., the polynomial $p$  must have degree at most three, and the space of negative weight deformations $Def^{<0}(\mathbb{C}^2/\mathbb{Z}_2\times \mathbb{C})$ is thus finite dimensional (equal to four in such case).

This phenomenon is general \cite{Odaka2024Bubbling}. In particular, we see that the canonical minimizing valuation is pointing out to this finite dimensional space $Def^{<0}$ of deformations (of the possibly jumping tangent cone) which, thanks to \cite{SSBubbling, sun2023no} and works of Conlon-Hein (e.g., \cite{conlon2025asymptotics}) indeed should parametrize the bubbles. More precisely, it is then natural to expect the following picture, generalizing to the full bubble tree the Li's normalized volume algebraic computations of the tangent cones.

\begin{conjpic}[\cite{SSBubbling, de2024some, Odaka2024Bubbling}] \label{CP1}

Given a polarized family of (KE) K-polystable varieties $\mathcal{X}\rightarrow \Delta$ with klt singularities and a section $\sigma$, there exist two canonical (depending only on the family and section), and AG computable, affine pointed varieties $$(\bar\sigma \in \textbf{W}^{\sigma}) \leadsto W_{\sigma(0)} \; \mbox{and} \;  (\bar\sigma \in \textbf{B}^{\sigma}) \leadsto C_{\sigma(0)}(Y),$$ negative weight deformations of the weighted cone $W_{\sigma(0)}$ and tangent cone $C_{\sigma(0)}(Y)$ at $\sigma(0)\in X_0$. In particular:\\

\begin{itemize}
\item Such $\textbf{B}^{\sigma}$ is characterized by admitting an AC CY metric (with tangent cone at infinity equal to $C_{\sigma(0)}(Y)$) which is the first such pGH limit of KE metrics on $X_t$ along $\sigma$ (\textit{minimal bubble}).\\
\item Replacing (after base change) the central fiber locally with $\textbf{B}^\sigma$ and iterating a finite number of times, we get a collection of $(\textbf{W}^\sigma_i, \textbf{B}^\sigma_i)$ computing all bubbling $pGH$ limits along $\sigma$.\\

\item Varying the section $\sigma$, one can determine the full bubble trees $\tau_0$, whose combinatorial structure embedds in the Berkovich's analytification of the family $$[\tau_0] \hookrightarrow \mathcal{X}^{AN}_{\C((t))},$$ via the scaling valuations (as in the above examples).

\end{itemize}

\end{conjpic}

Heuristically, the above picture should provide a way of obtaining a canonical geometric expansion of the degenerating  stable algebraic varieties $X_t$ to $X_0$, near a section  $\sigma$: $$\mathcal{X}/_{\mathbb{C}((t))}=X_0+B^\sigma_{min} t^{\lambda_{min}}+ \dots  + B^\sigma_{deep} t^{\lambda_{deep}},$$  or, more generally, without fixing the section,

$$ \mathcal{X}/_{\mathbb{C}((t))}= X_0+\sum_i  \tau_i t^{\lambda_i},$$
where $\tau_{i}$ are trees, inductively (but finitely) defined as $$\tau_i=B_i+\sum_j B_{ij} t^{\lambda_{ij}},$$ etc...,  with $B_{ijk...}$ asympotitically conical bubbling AC Calabi-Yau spaces and indexes running through stratifications of the spaces. \\

It is worth pointing out that, due to the possible jumping of tangent cones, the bubbling spaces $B_i^\sigma$ above are not, in general, parametrized by the versal deformation space $Def(X_0)$ of $X_0$. Indeed, even in the very special cases when there are no-local-to-global obstructions in the deformation theory of a singular space, i.e.,  $X_0$
$$ Def(X_0)=Def^{es}(X_0)\times Def(\mathcal{V}_p),$$
where $Def^{es}(X_0)$ parametrizes equisingular deformations and $Def(\mathcal{V}_p)$ the versal space of deformation of singularities, we see that in order to describe the  canonical geometry of the KE/K-stable space $X_t$ near $X_0$ it is more natural to consider $Def^{<0}(C(Y))$ instead of $Def(\mathcal{V}_p)$ (as moduli of first minimal bubbles), and so on. In analogy with the flat metric example, these considerations point toward the need of investigations and studies of possible modifications of the moduli spaces of KE/K-stable varieties which should parametrize bubblings (i.e, to the existence of certain notions of non-collapsing \textit{multiscale} K-moduli spaces).\\

From a more differential-geometric perspective, this conjectural picture leads to the possibility of obtaining computability of bubblings via Algebraic Geometry, via some sort of ``multiscale moduli continuity method'', that is, by realizing Sun's algebraic bubble of a given family via a-priori detectable algebraic rescalings. This can remove the need of performing  analytic gluing constructions to see (at least certain aspects) of bubbling. 

As an example of this, thanks to the already existent theory, we can see that the Stenzel's AC metric on the smoothing of the ODP point $\sum_i z_i^2=t$ must appear in degeneration of compact KE manifolds to nodal KE varieties. Indeed,  first one can compute the tangent cone algebraically, which gives the Stenzel's cone metric on the ODP. Sun's result gives a bubble with such asymptotic, which by Conlon-Hein theory must be a deformation of ODP. However, the versal deformation is just one dimensional, so the analytic structure of the bubble must the one of  $\sum_i z_i^2=t$, and the metric must be the Stenzel's one by uniqueness. So we have the following:

\begin{observation}[\cite{de2024some}] If a non-collapsed limit $X_0$ have only ODPs, then the Stenzel's AC CY metric on $T^*S^n$ must appear as (minimal equal deepest) bubble.
\end{observation}

If the conjectural picture holds, similar arguments will (re)construct Y.Li \cite{li2019new}, Szekelyhidi \cite{S19}, Conlon-Rochon \cite{conlon2025asymptotics} asymptotically conical CY type metrics with singular tangent cones at infinity, directly as bubbles from degenerating families of \textit{compact} KE manifolds.\\

In any case, it would be still interesting to investigate how this concrete algebraic bubbling  computability relates to the previously described deformation $Def$  theory, in order to possibly establish a deeper Kuranishi's theory for KE/K-moduli spaces at the boundary of the moduli in relation to more universal gluing constructions and to the geometric regularity of the degenerating KE metrics.\\

Recently Odaka \cite{Odaka2024Bubbling} has made some first (local!) advances for establishing the conjectural picture. For instance, he gave an abstract algebraic way to compute, possibly depending on certain choices, spaces $W_\sigma$ and $B_\sigma$, using a generalization of Halpern-Leistner $\Theta$-stratification induced by $\hat{vol}$ on the stack of pointed germs of klt singularities, which should compute the bubbling of a degenerating family of (compact) polarized KE manifolds. Indeed, his picture matches our rescalings for the two dimensional $A_k$ and flat $1$-dim cases. However, the full understanding of how this abstract local algebraic description relates to the global geometric one is still missing.

\section{Collapsing of polarized KE varieties:  moduli spaces compactifications and singularities} \label{S3}

This last section focuses on aspects of \textit{collapsing} of KE metrics in relation to Algebraic Geometry, and it consists of two (related) parts. In the first one, we aim to discuss the case of collapsing of polarized Ricci-flat Calabi-Yau manifolds. In particular, after reviewing the basically understood cases of limits of KE metrics for positive or negative KE metrics in relation to construction of compact moduli spaces of varieties, we will discuss and propose a  possible \textit{unifying picture} of compact moduli of polarized algebraic varieties when the canonical class has a sign (i.e., $K_X<,=,> 0$), in relation to metric limits at \textit{generic volume scale} of such canonical KE metrics. In the second part, instead, we will focus on the description of the collapsing part of negative KE metrics, discussing a conjectural picture inspired by the collapsing geometry of Calabi-Yaus.

\subsection{Unifying aspects of moduli compactifications for polarized KE varieties}

Compact algebraic moduli spaces of (smoothable) $\pm K_X$ polarized varieties $\overline{\mathcal{M}}^K$, parametrizing K-polystable varieties, have been constructed in recent years. \\

\subsubsection{Moduli of Fano varieties}
The Fano case (that is, $-K_X$ ample) is strongly related to the Yau-Tian-Donaldson conjecture \cite{chen_donaldson_sun2015_III}. In particular, it is known that there exist (abstract) projective algebraic compact moduli spaces of K-polystable Fano varieties, which are homeomorphic to the  Gromov-Hausdorff compactification (at volume, equal to diameter, scale) of the moduli spaces of smooth K\"ahler-Einstein/K-polystable Fano manifolds, always not collapsing: 
$$\textbf{ $\overline{\mathcal{M}}^K \cong \overline{\mathcal{M}_{KE>0}}^{GH}$}.$$ 

While the original constructions used the interplay between the algebraic geometry and the analytic properties the KE metrics (see, for instance, \cite{odaka2016compact, SpottiSunYao2016, Odaka2015CompactModuli, LiWangXu2019}), more recently a purely algebraic construction of $\overline{\mathcal{M}}^K$ has been obtained, including the non-smoothable situation, by Xu and collaborators in a series of papers (\cite{LiuXuZhuang2022} for the final one), and many examples, for instance in dimension $3$  (e.g., \cite{AraujoCastravetCheltsovEtAl2023}), are now becoming concretely computable. 

Furthermore, the existence of weak KE metrics on any singular K-polystable Fano variety has been shown to hold by the work of Berman, Boucksom and Jonsson using pluripotential variational techniques \cite{berman2021variational}. With this approach, however, the higher regularity of the weak KE metrics, which is needed to formulate statements on the metric Gromov-Hausdorff convergence as in the smoothable situation (for instance, if such spaces are RCD, e.g., \cite{szekelyhidi2024singular, fu2025rcd}), is currently  still missing.

\subsubsection{Moduli of $K_X$-ample varieties} KE Fano manifolds do not collapse. Thus, the above results conclude the investigation of the behavior of the KE metrics at volume (= diameter) scale. However, the case of ample canonical bundle $K_X>0$ is different. 

Beside non-collapsing situations similar to the Fano case (resulting in a \textit{partial} compactification of the moduli spaces adding KE varieties with klt singularities), now it is possible to have local collapsing at fixed volume scale, as the famous case of cusp formations of hyperbolic curves shows. Combining the works \cite{Kollar2013Dollar, Odaka2012CalabiKStability, Odaka2013, BermanGuenancia2014, SongSturmWang2020} we have that, there exists a projective algebraic  compactifications made of semi-log-canonical varieties (equivalently, K-polystable)  with ample canonical class, which corresponds to varieties admitting (weak) KE metrics:

$$\overline{\mathcal{M}}^K=_{Odaka}\overline{\mathcal{M}}^{KSBA} \cong\overline{\mathcal{M}_{KE<0}}^{GH, vol}.$$

In particular, \cite{SongSturmWang2020} shows that such identification happens at the volume scale in a \textit{generic} (that is, away from the volume small local collapsing regions) pointed Gromov-Hausdorff sense. Study of the local collapsing regions as been done by Zhang in \cite{ZhangKE}, and (tropical) aspects of the diameter scale limits by Odaka for the curves' case \cite{Odaka2014}.

\subsubsection{Moduli of polarized Calabi-Yau varieties}
The above Fano and $K_X$-ample cases (which also extend, for the most part, to \textit{pairs} $(X,D)$) provide an essentially clear picture of the relation between algebro-geometric compactifications of moduli spaces and degenerations of KE metrics (at the volume scale!). But what does it happen when we consider KE metrics that are Ricci-flat, i.e., in the Calabi-Yau case? That is:

\begin{problem}
    Does a canonical (algebraic!) compactification for polarized Calabi-Yaus (Ricci flat) varieties exist which encodes information about the degenerations of the canonical Ricci-flat metrics? 
\end{problem}

Compared to the previous cases, the fundamental new difficulty is that full collapsing to lower dimensional spaces can happen (indeed, consider the famous SYZ description of Mirror symmetry near maximal degenerations \cite{StromingerYauZaslow1996}). In the next, I want to describe some examples that, possibly, point towards the existence of a  unifying picture for moduli of polarized KE varieties which indeed includes such Calabi-Yau case as well, and it focuses on the (generic) fixed volume scale GH limits. However, before discussing it, it is quite important to recall the better (but not fully) understood partial compactifications and tropical limits (at the diameter scale).\\

First of all, we can have \textit{non-collapsing} limits of such Calabi-Yau metrics:  this is completely analogous to the KE case  above,  and gives a partial compactification of the moduli consisiting  of normal polarized Calabi-Yau varieties with canonical singularities (with singular KE metrics, in particular in the sense of \cite{EGZ2009}), again via Donaldson-Sun theory (see, e.g., \cite{ZhangCY}).  Moreover, such non-collapsed varieties are found precisely at finite WP distance by a result of Tosatti \cite{Tosatti2015}.\\ 

Some examples:\\

\textit{1. Flat 1D}. Our initial example of (conical) flat metric on $\mathbb{P}^1$: $$\mathcal{M(\underline{\beta})}\subsetneq \overline{\mathcal{M(\underline{\beta})}}^{Non-collaps} \cong \overline{\mathcal{M}}^{GIT, stable}_{\underline{\beta}} (\subsetneq  \overline{\mathcal{M}}^{GIT}_{\underline{\beta}}).$$
 here the non-collapsing  boundary is made by making points to come together, and it agrees with a GIT stable strata (collapsing polystable cusps are missing).\\
 
\textit{2. Polarized Abelian varieties}. In this special situation there is no non-smooth non-collapsing. Thus the partial compactification  agrees with the moduli of smooth (principal) polarized abelian varieties itself, and it is a locally symmetric space, via identification using periods  $$\mathcal{A}=\overline{\mathcal{A}}^{Non-collaps} \cong \Omega_A/\Gamma.$$\\ 

\textit{3. Polarized K3s}. As we discussed previously, Calabi-Yau metrics on K3 can degenerate in a non-collapsing way forming orbifold singularities. Thus we obtain a partial compactification adding such polarized K3 with ADE singularities (e.g., \cite{anderson2010survey}). Remarkably, after such ADE locus is added, we have description of the moduli partial compactification via periods as a locally symmetric space again:  $$\mathcal{F}_{2d} \subsetneq \overline{\mathcal{F}_{2d}}^{Non-collaps} \cong \Omega_F/\Gamma.$$ \\

For general situations (e.g., Calabi-Yau 3-folds) such partial compactification is obtained by adding normal Calabi-Yau varieties with canonical singularities. However, concrete examples parametrizing \textit{all} non-collapsing degenerations of given Calabi-Yau families seem still unknown.\\

The above concludes the description of the non-collapsing part. But what does it happen when collapsing do occur? In such situation the role of scale is important. The most understood and studied case is the one of \textit{fixed diameter} scale. For instance, this is the scale in which the SYZ description -of collapsing special Lagrangian torus fibrations at the large complex structure (LCS) limits- manifests \cite{StromingerYauZaslow1996}. First of all, at such scale we immediately have metric pre-compactness, by Gromov's theorem, so it make sense to talk about $\overline{\mathcal{M}}^{GH, diam}$.

However, metric limits at such scale won't be complex algebraic objects at all (indeed they will be real tropical varieties). Slightly more precisely, the guiding conjecture in the area (which generalizes away from the LCS limits a famous conjecture of Kontsevich-Soibelman \cite{KontsevichSoibelman2006}) is the following:
  
 \begin{conjecture}[Folklore conjecture generalizing Konstevich-Soibelman]
  The collapsing boundary $\partial_{collaps}\overline{\mathcal{M}}^{GH, diam}$ densely parametrizes singular tropical metric spaces homeomorphic to essential skeleta  $Sk(\mathcal{X})$ of algebraic degenerations (e.g., dual complexes of dlt degenerations). The generic boundary points corresponds to LCS limits. The space  $\overline{\mathcal{M}}^{GH, diam}$ is far from an algebraic object (indeed, the collapsing boundary $\partial_{collaps}\overline{\mathcal{M}}^{GH, diam}$ is itself a tropical object).
 \end{conjecture}

There is a very large literature on the topic, especially  at LCS limits (equivalently, when $dim(Sk(\mathcal{X}))=n$), for instance works of Gross-Seibert, Gross-Wilson, Tosatti, Zhang, etc... (e.g.,  \cite{GrossSiebert2011, GrossWilson2000}). More recently, there have been works also for $dim(Sk(\mathcal{X}))\leq n$  by e.g.,  Y. Li \cite{Li2022, Li2023},  Sun-Zhang \cite{SunZhang2019}, Boucksom-Jonsson \cite{BoucksomJonsson2017}, etc... 

However, what can we say regarding such tropical compactifications in our three main examples above? Mostly thank to the works of Odaka, the answer is completely known in such cases.\\

\textit{1. Flat 1D}. In such case the tropical GH compactification (not trivial only for non-generic value of cone angles $\underline{\beta}$s) consists in adding \textit{points} at the cusps, parametrising unit \textit{intervals}: $$\overline{\mathcal{M(\underline{\beta})}}^{GH,diam} \cong_{top} \overline{\mathcal{M}}^{GIT}_{\underline{\beta}}.$$ 
Note that although  $\overline{\mathcal{M}}^{GIT}_{\underline{\beta}}$ happens to be algebraic, its identification at the cusp is purely topological, that is, at such scale there is no any natural universal modular algebraic structure. \\

\textit{2. Polarized abelian varieties}. For polarized abelian varieties one adds real lower dimensional tori. It turned out (\cite{Odaka2019}) that the space parametrizing such limits is given by a purely topological Stake-type compactification (called Satake adjoint) coming from the representation theory of the locally symmetric space description $$\overline{\mathcal{A}}^{GH, diam} \cong  \overline{\Omega_A/\Gamma}^{Sat-Adj}.$$  \\

\textit{3. Polarized K3s}. Remarkably, thanks to work of Odaka-Oshima \cite{OdakaOshima2021}, the same holds true for K3 by adding intervals and tropical K3s (certain affine structures on the sphere): $$\overline{\mathcal{F}_{2d}}^{GH, diam} \cong  \overline{\Omega_F/\Gamma}^{Sat-Adj}.$$  More recentely, there have been further refinments of the above parametrizing \textit{measure}-GH compactifications (see, for instance, \cite{Odaka2022a}, \cite{HondaSunZhang2019}).\\

Beside the above, there are not yet other known examples of complete descriptions of moduli spaces at such diameter scale. Moreover,  beyond  some tentative conjectural picture in \cite{OdakaOshima2021}, there is not yet a clear general formulation of what the tropical boundary of the moduli, parametrizing tropical varieties (GH limits), would be in relation to algebraic geometry. \\

Finally, we have reached the next universally relevant scale, that is the fixed volume scale (which is indeed quite natural in relation to fixing the polarization). In our three examples the situation looks as follow:
\\

\textit{1. Flat 1D}. By adding polystable point to parametrize \textit{lines} we get the algebraic compactification: $$\overline{\mathcal{M(\underline{\beta})}}^{GH,vol} \cong \overline{\mathcal{M}}^{GIT}_{\underline{\beta}}$$ \\

\textit{2. Pol. Abelian}. In \cite{Odaka2019}, Odaka showed that limits at constant volume scale are of the form $\mathbb{R}^k\times E$, $E$ lower dimensional abelian variety and (using the explicit description in term of locally symmetric spaces): $$\overline{\mathcal{A}}^{GH, vol} \cong  \overline{\Omega_A/\Gamma}^{BB},$$ where BB denotes the classical algebraic Baily-Borel compactification. Let's see this in a trivial example. Just take the product of a fixed elliptic curve $E_0$ with the universal one $E_\tau$.  If we look to the metric limit at fixed volume scale (after fixing a polarization) of the flat metrics on $X_t=E_0 \times E_\tau$ it is clear that the metrics on the $E_\tau$ factor stretch converging to a real line while the flat metric on $X_0$ remains untouched. Thus the volume scale limit is just $E_0\times \R$ (which indeed also matches the algebro-geometric meaning of the Baily-Borel compactification, whose boundary in such case is just the universal modular curve parametrising the splitting compact ellitic curves factors arising in semi-abelian degenerations). \\

\textit{3.} What does it happen for polarized K3 surfaces? I claim that also in this case the natural volume scale compactification is indeed still the Baily-Borel (for the locally symmetric space $\Omega_F/\Gamma$). From some explicit hypersurface examples $\mathcal{X}\rightarrow \Delta$, with $dim(Sk(\mathcal{X}))=1$, for instance via gluings  \cite{HeinSunViaclovskyZhang2020}, \cite{SunZhang2019}), we can see that at volume scale the geometry looks \textit{generically} like $ E\times \mathbb{R}$, with $E$ flat elliptic curve $E= X_0^1 \cap  X_0^2 $ \textit{lc center or source} of a simple normal crossing central fiber. More generally, as an application of the  classification of collapsing K3s via Gibbons-Hawking ansatz by Sun-Zhang \cite{SunZhang2024} combined with the understanding of the Baily-Borel compactification in terms of periods (see, e.g., Friedman-Scattone \cite{Friedman1986}), we see the following:

\begin{observation} The boundary points of the Baily-Borel compactification $\overline{\Omega_F/\Gamma}^{BB}$ parametrize the generic pGH limits at the volume scale for Ricci-flat metrics on polarized K3s:
\begin{itemize}
\item $E\times \mathbb{R}$, with $E$ appearing as lc center of type II models ($1$-dim BB strata); 
\item $\mathbb{R}^2$, for type III models (isolated cusps in BB).
\end{itemize}
\end{observation}

Here a sketch of the proof. Using Sun-Zhang results in combination with periods computations (which are possible using their a-priori ansatz description of the metric) for a degenerating family one shows that in the polarized case there cannot be collapsing to real three dimensional tori (e.g., as the ones appearing in  Foscolo's gluing \cite{Foscolo2019}). If the limit at diameter scale is two dimensional, then periods show that we necessarily limiting at the cusps of the BB compactification and that, generically, the metric looks like a shrinking torus bundle  (that is, we are at a smooth general fiber of a the special Lagrangian fibration). Taking  now the volume scale, we see that the torus fibration still shrinks while the base expands, thus having a flat $\mathbb{R}^2$ as generic volume limit (one should think to it as $\{p\}\times \mathbb{R}^2$, with the point $p$ identified as a log canonical center of a type III filling, e.g., the point of intersection of three irreducible components in a snc degeneration).

If the limit is instead an interval, one rule out cylindrical splittings via periods as well, and thus only Heisenberg symmetries survive. Using the Gibbons-Haking's ansatz description of the metric in such case, one shows that at volume scale the metric looks like generically as a flat metric on a torus times the real line. Finally, using results of Friedman-Scattone on the geometric meaning of the Baily-Borel compactification in relation to periods and Kulikov models, one can identify such flat torus with the log canonical center of a type II degeneration. This observation can be viewed as an extension to the BB boundary of the metric proof of surjectivity of the period map by Liu \cite{Liu2023}.

\begin{remark} Note that in the above compactification we are adding a plane to each cusps (and not identifying them all together). Similarly, some flat elliptic curve appears as compact splitting factor  more than one time at the boundary (the one dimensional components at the boundary are indeed covers of the moduli of elliptic curves). 

Heuristically, the BB compactification is the normalization of the generic volume scale GH compactification.  Observe also that  in all such three examples of compactifications at volume scale the moduli spaces are \textit{projective algebraic}.
\end{remark}

Optimistically, this leads to formulate the following tentative conjectural picture, which would definitely require some fine tuning.

\begin{conjpic} \label{CP2}

For polarized families $\mathcal{X}^*\rightarrow \Delta^*$ of $n$-dim CYs, the \textbf{\textit{generic} pGH limits at volume scale} of polarized Ricci-flat metrics are of the form:
$$Y^{n-s}\times \mathbb{R}^{s},$$
with $s=dim(Sk(\mathcal{X}))$, and $Y^{n-s}$ a $(n-s)$-dim polarized CY with canonical singularities, appearing as minimal \textit{lc center} of (non-unique) minimal dlt models (compact CY at the intersections of components $\bigcap_{i=1}^{s+1} X_0^i$).

\end{conjpic}

Several remarks are needed. First of all, the precise meaning of \textit{generic} need to be clarified, possibly in relation to limits of natural measures as studied in the works of Boucksom-Jonsson \cite{BoucksomJonsson2017}. Similarly, some care need to be taken in relation to limiting polarizations and flopping ambiguity of the lc center. Beside the above K3 and abelian varieties examples, this conjecture seems compatible with Y. Li, Sun-Zhang and Odaka's works previously mentioned. In particular, such notion of \textit{generic} pGH limit at volume scale provide a tentative refinement/answer to  Question B.8  in \cite{Odaka2022b}. Obtaining higher estimates in Y. Li's works would show the validity of the conjecture for the families he considers.\\

In relation to canonical moduli compactifications, this suggests the existence of a compactification parametrizing compact splitting components of generic Gromov-Hausdorff limits at volume scale, which algebro-geometrically identify with log canonical centers of degenerations. 

The boundary components of the moduli would themselves be moduli of certain lower dimensional Calabi-Yau varieties:

$$\overline{\mathcal{M}}^{GH,vol} \cong\overline{\mathcal{M}}^{non-collaps} \bigsqcup \{\mbox{(covering of) moduli of lower dim pol. CYs}\}=: \overline{\mathcal{M}}^{K}.$$

Such moduli compactifications  $\overline{\mathcal{M}}^{K}$ should be \textit{projective algebraic} varieties.  The projectivity could follow by study the extension of the CM bundle, also in relation to positivity and limiting aspects of the Weil-Petersson metric.

Morally, this could be seen as a generalized Baily-Borel type compactification (carrying a modular meaning in relation to metric collapsing of the Ricci-flat metrics at volume scale). The very recent  work by Bakker-Filipazzi-Mauri-Tsimerman \cite{BakkerFilipazziMauriTsimerman2025} seems to provide a quite similar Hodge-theoretical compactification, and it would then be nice to explore its possible connections with the differential geometric picture above.\\

This notion of (generic) limit at volume scale would then \textit{unify }the notion of canonical algebraic moduli compactifications of moduli of polarized varieties in relation to limits of KE metrics, so giving a natural extension of what a K-moduli could be also in the presence of full collapsing. In this regard, it could be very interesting to analyze what happens to the above picture for the pair case $(X,\beta D)$, since such K-moduli space of Calabi-Yaus  would fit between the two better understood K-moduli of Fano and $K_X$-ample cases and thus it should be the natural limit of such moduli. For instance, a first example could be to study from this point of view the Baily-Borel related moduli at the CY threshold for surfaces pairs considered in, e.g., \cite{AscherEtAl2023}.

\begin{remark}\label{HS}
    What does it happen at even higher scales for moduli of Calabi-Yaus? For instance, what if we rescale the metric just enough to avoid any collapsing (although, specify such scale \textit{universally} in the moduli, contrary to the above fixed volume case, may not be possible)? This will lead us to some interesting limits given by complete Calabi-Yau varieties collapsing at infinity (collapsing bubbles). As examples suggest (e.g., \cite{SunZhang2019}), we will see then the emergence of Tian-Yau metrics and its generalization \cite{CollinsLi2024}. Such spaces can be compactified as irreducible components of certain degenerations, however, in general, their relation with algebraic geometry remains unclear, beside some interesting preliminary works by Odaka \cite{Odaka2022b}, and a full picture missing even in the case of polarized K3 (despite, thanks to the works of Sun and Zhang \cite{SunZhang2024}, we have now a very good a-priori understanding of the metric limits as (classified) gravitational instantons). 
    
    A final remark/warning: separating such minimal non-collapsing scales and the AC conical bubbling pictures we describe in the previous section may not be possible. For instance, consider a product of a family of non-collapsing CY (where AC bubbling occurs) with a family of collapsing one (where minimal Tian-Yau bubbling occurs): at the scale when see the emergence of Tian-Yau spaces, we will also see the emergence of conical bubbles from the first factor, thus tangling the two pictures together generating interesting non-collapsed bubbles' trees worth to be studied algebro-geometrically as well.

\end{remark}

\subsection{Algebro-geometric aspects of isolated log-canonical KE cusps}

We now move to the final part of this survey, where we discuss an example of collapsing near singularities of negative KE metrics. Even though at first one may think this has not much to do with the previous Calabi-Yau collapsing picture, I want to describe some concrete examples that suggest instead an analogy between such two cases in the description of the collapsing geometries. The general problem is the following:

\begin{problem}
Understanding  (multiscale) collapsing behaviors of negative KE near isolated log canonical cusps in relation to Algebraic Geometry.
\end{problem}

Beside the general results of Berman-Guenancia on the existence of negative KE metrics on semi-log-canonical varieties \cite{BermanGuenancia2014}, the geometric understanding of the metric near singularities (away from the non-collapsing klt case) in connection to the underlying algebraic geometry remains quite weak, even conjecturally. Special cases of collapsing are known, mostly from a more differential geometric viewpoint, thanks to some recent works (Datar-Fu-Song \cite{DatarFuSong2023},  Fu-Hein-Jiang \cite{FuHeinJiang2021}, etc.), but a general \textit{framework} to describe what happens seems missing.\\

Let us review some concrete examples, which connect to the above works.\\

\textit{1. Algebraic cones.} These very special lc singularities (the contraction of the zero section of the total space of the dual of an ample line bundle on a Calabi-Yau manifold $X$) can be canonically \textit{metrized} via a Calabi's ansatz:
 $$g=\frac{1}{2}dt^2 + e^{-t}(\pi^*g_{CY} + 2e^{-t}d\theta^2)\rightarrow_{t\rightarrow \infty} \mathbb{R}.$$

The WP cuspidal regions in $\overline{\mathcal{M}(\underline{\beta})}^{GH}$ we described in the Section \ref{S1} are indeed an example of such singularities (for pairs). In particular, in full analogy with such example, horospheres $(t=c)$ fully collapse as well, and the CY metric on $X$ (in the Kähler class of the polarizing line bundle) emerges at higher scale. However, this is a purely local picture. Thus, what happens when such singularity is found in a compact singular KE variety? The expectation is that any KE metric is indeed asymptotic to such local model. This has been recently confirmed to be true in the case when $X$ is a flat torus \cite{FuHeinJiang2021}.\\

\textit{2. Hilbert modular cusps, 2d lc}. This example is given by arithmetic quotients  $(\mathbb{H}\times\mathbb{H})/\Gamma$, and the KE metric is simply the one induced by the product of the two hyperbolic factors. A concrete analysis shows that such negative KE metrics collapse at infinity to a cylinder $S^1\times \mathbb{R}$. It has been recently shown that any complete negative KE metric near such singularity is close to such model \cite{DatarFuSong2023}.\\

 The above two cases essentially conclude the two dimensional case (see also the picture suggested by Kobayashi \cite{Kobayashi1985}). However, what do these examples have in common from an algebro-geometric point of view? Giving an answer to such question is relevant in order to formulate some tentative general picture to study higher dimensions as well. Note that that in the first example the total space of the line bundle provides the minimal resolution of singularity, with the Calabi-Yau manifold $X$ as the exceptional divisor. In the second case, the minimal resolution of such modular cusps is instead given by an exceptional cycle of rational curves,  as famously described by Hirzebruch. If we consider the dual complexes of these two minimal resolutions, we get a single point in the first case and a circle  $S^1$ in the second one. Thus, the following result by Engberg, based on explicit studies, gives then a possible answer to our question.

\begin{theorem}[\cite{Engberg2022}] Horospheres (for the model KE metrics on modular cusps) naturally collapse to the dual complex of the minimal resolution.
\end{theorem}

Here naturally means by composing with a Boucksom-Jonsson log map (\cite{BoucksomJonsson2017}). For instance, one can see that the $S^1$-fibration on the link of the singularity naturally focuses (as we get closer to the singularity) at the intersection of the rational curves. In particular, the $S^1$ factor in the asymptotic cylinder can be canonically identified with the dual complex of the minimal resolution, hence giving an algebro-geometric meaning to it. \\

We can now try to use such two dimensional picture to develop some higher dimensional framework. Firstly, let's subdivide isolated lc singularities according to the dimension $d\leq{n-1}$ of dual complex $\Delta^{d}$ of a minimal dlt resolutions. This is well defined thanks to work by De Fernex, Koll\'ar and Xu \cite{deFernexKollarXu2017}. For the case $d=0$, as we said, we have an algebraic cone like picture, as in the Calabi's ansatz example given before. The generalization of the second example would instead be  given by the maximal $d=n-1$ dimensional case, basically the analogous of large complex structure limits for Calabi-Yaus. In particular, considering certain higher dimensional generalizations of the two dimensional case leads to the following result (also based on concrete computations), which confirms in higher dimensions Engberg's interpretation:

\begin{theorem}[\cite{FangSpotti25}] For $n$-dim Tsuchihashi cusps with model KE metrics (Cheng-Yau), we have the following collapsing picture at infinity:

\begin{itemize}
\item Horospheres collapse to the dual complex $\Delta^{n-1}$ of a natural toroidal resolution;
\item The KE metrics collapse to metrics on $\Delta^{n-1}\times \mathbb{R}$, solutions of a real Monge-Amp\`ere equation $det(\partial^2v)=e^{(n+2)v}$.
\end{itemize}

\end{theorem}

Optimistically, one is then tempted to expect a picture for describing the collapsing of negative KE cuspidal metric similar to the Calabi-Yau degenerations for each dimension $d$ of the dual complex $\Delta^d$:

\begin{conjpic}\label{CP3} Complete negative KE metrics near isolated lc cusps collapse to tropical metrics on $\Delta^d\times \mathbb{R}$, as canonical limitig data.
\end{conjpic}

This can be viewed as a local analogous of the (generalized) Kontsevich-Soibelman conjecture, and it would answer a question raised in Datar-Fu-Song \cite{DatarFuSong2023} about the possible canonicity of the collapsing limits, playing a role somewhat similar to tangent cones at singularities in the non-collapsing case. \\

Continuing, more speculatively, with such Calabi-Yau analogy, we could hope to see the emergence at higher scales of splitting Calabi-Yau metrics on compact lc centers and, zooming even more, splitting complete CY metrics compactified as some exceptional components. In particular, in order to support  and refine this even more general picture, it would be interesting to develop local Y. Li and Sun-Zhang's type examples, for instance in complex dimension three.\\

In conclusion, in this survey I presented some old and new examples  describing how Algebraic Geometry could be used to understand degenerations of KE metrics. As I hope it is clear to the reader, these types of studies should provide newer point-of-views in the theory of moduli spaces compactifications and singularities' formations, areas which are still full of fundamental questions that a \textit{combined} differential and algebraic viewpoint can help to unfold.

\newpage

\bibliographystyle{alpha}

\bibliography{MultiKE}

\newcommand{\etalchar}[1]{$^{#1}$}
\begin{thebibliography}{BFMT25}

\bibitem[ABB{\etalchar{+}}23]{AscherEtAl2023}
Kenneth Ascher, Dori Bejleri, Harold Blum, Kristin DeVleming, Giovanni Inchiostro, Yuchen Liu, and Xiaowei Wang.
\newblock {Moduli of boundary polarized Calabi--Yau pairs}.
\newblock {\em arXiv preprint arXiv:2307.06522}, 2023.

\bibitem[ACC{\etalchar{+}}23]{AraujoCastravetCheltsovEtAl2023}
Carolina Araujo, Ana-Maria Castravet, Ivan Cheltsov, Kento Fujita, Anne-Sophie Kaloghiros, Jesus Martinez-Garcia, Constantin Shramov, Hendrik Süß, and Nivedita Viswanathan.
\newblock {\em {The Calabi Problem for Fano Threefolds}}, volume 485 of {\em London Mathematical Society Lecture Note Series}.
\newblock Cambridge University Press, 2023.

\bibitem[And89]{Anderson1989}
Michael~T. Anderson.
\newblock {Ricci curvature bounds and {E}instein metrics on compact manifolds}.
\newblock {\em Journal of the American Mathematical Society}, 2(3):455--490, 1989.

\bibitem[And10]{anderson2010survey}
Michael~T. Anderson.
\newblock {A Survey of Einstein Metrics on 4-Manifolds}.
\newblock In Shing-Tung Yau, editor, {\em Handbook of Geometric Analysis, Volume 3}, pages 1--102. International Press, 2010.

\bibitem[BBJ21]{berman2021variational}
Robert~J. Berman, Sébastien Boucksom, and Mattias Jonsson.
\newblock {A variational approach to the Yau–Tian–Donaldson conjecture}.
\newblock {\em Journal of the American Mathematical Society}, 34(3):605--652, 2021.

\bibitem[BFMT25]{BakkerFilipazziMauriTsimerman2025}
Benjamin Bakker, Stefano Filipazzi, Mirko Mauri, and Jacob Tsimerman.
\newblock {Baily--Borel compactifications of period images and the b-semiampleness conjecture}.
\newblock arXiv preprint arXiv:2508.19215, 2025.

\bibitem[BG14]{BermanGuenancia2014}
Robert~J. Berman and Henri Guenancia.
\newblock {K\"ahler–Einstein metrics on stable varieties and log canonical pairs}.
\newblock {\em Geometric and Functional Analysis}, 24(6):1683--1730, 2014.

\bibitem[BJ17]{BoucksomJonsson2017}
Sébastien Boucksom and Mattias Jonsson.
\newblock {Tropical and non-Archimedean limits of degenerating families of volume forms}.
\newblock {\em Journal de l’École polytechnique — Mathématiques}, 4:87--139, 2017.

\bibitem[BKN89]{BandoKasueNakajima1989}
Shigetoshi Bando, Katsumi Kasue, and Hiraku Nakajima.
\newblock {On a construction of coordinates at infinity on manifolds with fast curvature decay and maximal volume growth}.
\newblock {\em Inventiones Mathematicae}, 97:313--349, 1989.

\bibitem[Blu18]{Blum2018Existence}
Harold Blum.
\newblock {Existence of valuations with smallest normalized volume}.
\newblock {\em Compositio Mathematica}, 154(4):820--849, 2018.

\bibitem[BR14]{BiquardRollin2014}
Olivier Biquard and Yann Rollin.
\newblock {Smoothing singular extremal K{\"a}hler surfaces and minimal Lagrangians}.
\newblock {\em Annales de l'Institut Fourier}, 64(3):1127--1156, 2014.

\bibitem[CDS15]{chen_donaldson_sun2015_III}
Xiuxiong Chen, Simon Donaldson, and Song Sun.
\newblock {K\"ahler-Einstein metrics on Fano manifolds. III: Limits as cone angle approaches $2\pi$ and completion of the main proof}.
\newblock {\em Journal of the American Mathematical Society}, 28(1):235--278, 2015.

\bibitem[CH25]{conlon2025asymptotics}
Ronan~J. Conlon and Hans-Joachim Hein.
\newblock {Asymptotics for resolutions and smoothings of Calabi-Yau cones}.
\newblock {\em arXiv preprint arXiv:2503.16702}, 2025.

\bibitem[Che01]{cheeger_degen_riemannian_2001}
Jeff Cheeger.
\newblock {\em {Degeneration of Riemannian metrics under Ricci curvature bounds}}.
\newblock Publications of the Scuola Normale Superiore. Edizioni della Normale, Pisa, 2001.

\bibitem[CHL05]{CouwenbergHeckmanLooijenga2005}
Wim Couwenberg, Gert Heckman, and Eduard Looijenga.
\newblock {Geometric structures on the complement of a projective arrangement}.
\newblock {\em Publications Mathématiques de l'IHÉS}, 101:69--161, 2005.

\bibitem[CL24]{CollinsLi2024}
Tristan~C. Collins and Yang Li.
\newblock {Complete Calabi–Yau metrics in the complement of two divisors}.
\newblock {\em Duke Mathematical Journal}, 173(18):3559--3604, 2024.

\bibitem[dBS19]{deBorbonSpotti2019}
Martin de~Borbon and Cristiano Spotti.
\newblock {Local models for conical K{\"a}hler--Einstein metrics}.
\newblock {\em Proceedings of the American Mathematical Society}, 147(3):1217--1230, 2019.

\bibitem[dBS23]{deBorbonSpotti2023}
Martin de~Borbon and Cristiano Spotti.
\newblock {Calabi–Yau metrics with conical singularities along line arrangements}.
\newblock {\em Journal of Differential Geometry}, 123(2):195--239, 2023.

\bibitem[dBS24]{de2024some}
Martin de~Borbon and Cristiano Spotti.
\newblock {Some models for bubbling of (log) K{\"a}hler--Einstein metrics}.
\newblock {\em Annali dell'Universit\'a di Ferrara}, 70(3):1037--1068, 2024.

\bibitem[dFKX17]{deFernexKollarXu2017}
Tommaso de~Fernex, János Kollár, and Chenyang Xu.
\newblock {The dual complex of singularities}.
\newblock In K.~Oguiso, K.~Matsuki, and T.~Shioda, editors, {\em Higher Dimensional Algebraic Geometry—In Honour of Professor Yujiro Kawamata’s 60th Birthday}, volume~74 of {\em Adv. Stud. Pure Math.}, pages 103--130. Math. Soc. Japan, 2017.

\bibitem[DFS23]{DatarFuSong2023}
Ved Datar, Xin Fu, and Jian Song.
\newblock { K\"ahler–Einstein metric near an isolated log canonical singularity}.
\newblock {\em Journal für die reine und angewandte Mathematik}, 797:79--116, 2023.

\bibitem[DM86]{DeligneMostow1986}
Pierre Deligne and George~D. Mostow.
\newblock Monodromy of hypergeometric functions and non-lattice integral monodromy.
\newblock {\em Publications Mathématiques de l'IHÉS}, 63:5--89, 1986.

\bibitem[DS14]{donaldson_sun_gh_limits_2014}
Simon~K. Donaldson and Song Sun.
\newblock {Gromov--Hausdorff limits of K\"ahler manifolds and algebraic geometry}.
\newblock {\em Acta Mathematica}, 213(1):63--106, 2014.

\bibitem[DS17]{donaldson_sun_gh_limits_II_2017}
Simon~K. Donaldson and Song Sun.
\newblock {Gromov-Hausdorff limits of K\"ahler manifolds and algebraic geometry, II}.
\newblock {\em Journal of Differential Geometry}, 107(2):327--371, 2017.

\bibitem[EGZ09]{EGZ2009}
Philippe Eyssidieux, Vincent Guedj, and Ahmed Zeriahi.
\newblock {Singular K\"ahler–Einstein Metrics}.
\newblock {\em Journal of the American Mathematical Society}, 22(3):607--639, 2009.

\bibitem[Eng22]{Engberg2022}
Lukas Engberg.
\newblock {\em {On some algebraic aspects of special K\"ahler metrics}}.
\newblock PhD thesis, Aarhus University, 2022.

\bibitem[FGS25]{fu2025rcd}
Xin Fu, Bin Guo, and Jian Song.
\newblock {RCD structures on singular K\"ahler spaces of complex dimension three}.
\newblock {\em arXiv preprint arXiv:2503.08865}, 2025.

\bibitem[FHJ21]{FuHeinJiang2021}
Xin Fu, Hans‐Joachim Hein, and Xumin Jiang.
\newblock {Asymptotics of K\"ahler–Einstein metrics on complex hyperbolic cusps}.
\newblock {\em arXiv preprint arXiv:2108.13390}, 2021.

\bibitem[FHJ25]{fu2025continuous}
Xin Fu, Hans-Joachim Hein, and Xumin Jiang.
\newblock {A Continuous Cusp Closing Process for Negative K\"ahler–Einstein Metrics}.
\newblock {\em Geometric and Functional Analysis}, 35(2):542--632, mar 2025.

\bibitem[Fos19]{Foscolo2019}
Lorenzo Foscolo.
\newblock {ALF Gravitational Instantons and Collapsing Ricci-Flat Metrics on the K3 Surface}.
\newblock {\em Journal of Differential Geometry}, 112(1):79--120, 2019.

\bibitem[FS86]{Friedman1986}
Robert Friedman and Francesco Scattone.
\newblock {Type III Degenerations of K3 Surfaces}.
\newblock {\em Inventiones Mathematicae}, 83:1--40, 1986.

\bibitem[FS25]{FangSpotti25}
Yanbo Fang and Cristiano Spotti.
\newblock {Some examples of the collapsing behavior of K\"ahler-Einstein metrics near log-canonical singularities}.
\newblock {\em In preparation}, 2025.

\bibitem[GH78]{GibbonsHawking1978}
G.~W. Gibbons and S.~W. Hawking.
\newblock {Gravitational Multi-Instantons}.
\newblock {\em Physics Letters B}, 78(4):430--432, 1978.

\bibitem[GS11]{GrossSiebert2011}
Mark Gross and Bernd Siebert.
\newblock From real affine geometry to complex geometry.
\newblock {\em Annals of Mathematics}, 174(3):1301--1428, 2011.

\bibitem[GW00]{GrossWilson2000}
Mark Gross and P.~M.~H. Wilson.
\newblock {Large complex structure limits of K3 surfaces}.
\newblock {\em Journal of Differential Geometry}, 55(3):475--546, 2000.

\bibitem[Has03]{Hassett2003}
Brendan Hassett.
\newblock Moduli spaces of weighted pointed stable curves.
\newblock {\em Advances in Mathematics}, 173(2):316--352, 2003.

\bibitem[HSVZ20]{HeinSunViaclovskyZhang2020}
Hans-Joachim Hein, Song Sun, Jeff Viaclovsky, and Ruobing Zhang.
\newblock {Nilpotent Structures and Collapsing Ricci-Flat Metrics on K3 Surfaces}.
\newblock {\em Journal of the American Mathematical Society}, 33(3):755--823, 2020.

\bibitem[HSZ19]{HondaSunZhang2019}
Shouhei Honda, Song Sun, and Ruobing Zhang.
\newblock {A Note on the Collapsing Geometry of Hyperk{\"a}hler Four-Manifolds}.
\newblock {\em Science China Mathematics}, 62(11):2195--2210, 2019.

\bibitem[Kob85]{Kobayashi1985}
Ryoichi Kobayashi.
\newblock {Einstein–K\"ahler metrics on open algebraic surfaces of general type}.
\newblock {\em Tohoku Mathematical Journal, Second Series}, 37(1):43--77, 1985.

\bibitem[Kol13]{Kollar2013Dollar}
J{\'a}nos Koll{\'a}r.
\newblock Moduli of varieties of general type.
\newblock In {\em Handbook of Moduli. Vol. II}, volume~25 of {\em Adv. Lect. Math. (ALM)}, pages 131--157. Int. Press, Somerville, MA, 2013.

\bibitem[Kro89]{Kronheimer1989}
Peter~B. Kronheimer.
\newblock {The construction of ALE spaces as hyper-K{\"a}hler quotients}.
\newblock {\em Journal of Differential Geometry}, 29(3):665--683, 1989.

\bibitem[KS06]{KontsevichSoibelman2006}
Maxim Kontsevich and Yan Soibelman.
\newblock {Affine structures and non-Archimedean analytic spaces}.
\newblock In {\em The Unity of Mathematics}, volume 244 of {\em Progress in Mathematics}, pages 321--385. Birkhäuser, 2006.

\bibitem[Lan03]{Langer2003}
Adrian Langer.
\newblock {Logarithmic orbifold Euler numbers of surfaces with applications}.
\newblock {\em Proceedings of the London Mathematical Society (3)}, 86(2):358--396, 2003.

\bibitem[LeB15]{LeBrun2015}
Claude LeBrun.
\newblock {Einstein Metrics, Harmonic Forms, and Symplectic Four-Manifolds}.
\newblock {\em Annals of Global Analysis and Geometry}, 48(1):75--85, 2015.

\bibitem[Li18]{Li2018NormalizedVolumes}
Chi Li.
\newblock {Minimizing normalized volumes of valuations}.
\newblock {\em Mathematische Zeitschrift}, 289(1-2):491--513, 2018.

\bibitem[Li19]{li2019new}
Yang Li.
\newblock {A new complete Calabi–Yau metric on \( \mathbb{C}^3 \)}.
\newblock {\em Inventiones mathematicae}, 215(2):433--473, 2019.

\bibitem[Li21]{Li2021}
Chi Li.
\newblock {On the stability of extensions of tangent sheaves on K{\"a}hler--Einstein Fano/Calabi--Yau pairs}.
\newblock {\em Mathematische Annalen}, 381(3--4):1943--1977, 2021.

\bibitem[Li22]{Li2022}
Yang Li.
\newblock {Strominger–Yau–Zaslow conjecture for Calabi–Yau hypersurfaces in the Fermat family}.
\newblock {\em Acta Mathematica}, 229(1):1--53, 2022.

\bibitem[Li23]{Li2023}
Yang Li.
\newblock {Intermediate complex structure limit for Calabi–Yau metrics}.
\newblock {\em arXiv preprint arXiv:2305.02258}, 2023.

\bibitem[Liu23]{Liu2023}
Hongyi Liu.
\newblock {A short proof of the surjectivity of the period map on {K3} manifolds}.
\newblock {\em arXiv preprint arXiv:2304.09501}, 2023.

\bibitem[LM25]{LiMiao2025VolumeKSemistableFano}
Chi Li and Minghao Miao.
\newblock {On the volume of K-semistable Fano manifolds}.
\newblock {\em arXiv preprint arXiv:2508.19215}, 2025.

\bibitem[LWX19]{LiWangXu2019}
Chi Li, Xiaowei Wang, and Chenyang Xu.
\newblock {On the proper moduli spaces of smoothable K\"ahler–Einstein Fano varieties}.
\newblock {\em Duke Mathematical Journal}, 168(8):1387--1459, 2019.

\bibitem[LWX21]{LiWangXu2021TangentCones}
Chi Li, Xiaowei Wang, and Chenyang Xu.
\newblock {Algebraicity of metric tangent cones and equivariant K-stability}.
\newblock {\em Journal of the American Mathematical Society}, 34(4):1175--1214, 2021.

\bibitem[LX19]{liu2019k}
Yuchen Liu and Chenyang Xu.
\newblock K-stability of cubic threefolds.
\newblock {\em Duke Mathematical Journal}, 168(11):2029--2073, 2019.

\bibitem[LXZ22]{LiuXuZhuang2022}
Yuchen Liu, Chenyang Xu, and Ziquan Zhuang.
\newblock {Finite generation for valuations computing stability thresholds and applications to K-stability}.
\newblock {\em Annals of Mathematics}, 196(2):507--566, 2022.

\bibitem[McM17]{McMullen2017}
Curtis~T. McMullen.
\newblock {The Gauss–Bonnet theorem for cone manifolds and volumes of moduli spaces}.
\newblock {\em American Journal of Mathematics}, 139(1):261--291, 2017.

\bibitem[MM93]{MabuchiMukai1993}
Toshiki Mabuchi and Shigeru Mukai.
\newblock {Stability and Einstein--K\"ahler metric of a quartic del Pezzo surface}.
\newblock In {\em Einstein Metrics and Yang--Mills Connections}, volume 145 of {\em Lecture Notes in Pure and Applied Mathematics}, pages 133--160. Marcel Dekker, New York, 1993.

\bibitem[MSY06]{martelli2006geometric}
Dario Martelli, James Sparks, and Shing-Tung Yau.
\newblock {The geometric dual of a-maximisation for toric Sasaki–Einstein manifolds}.
\newblock {\em Communications in Mathematical Physics}, 268(1):39--65, 2006.

\bibitem[MZ20]{MazzeoZhu2020}
Rafe Mazzeo and Xuwen Zhu.
\newblock {Conical metrics on Riemann surfaces I: The compactified configuration space and regularity}.
\newblock {\em Geometry \& Topology}, 24(1):309--372, 2020.

\bibitem[NO25]{namikawa2025canonical}
Yoshinori Namikawa and Yuji Odaka.
\newblock {Canonical Torus Action on Symplectic Singularities}.
\newblock {\em arXiv preprint arXiv:2503.15791}, 2025.

\bibitem[Oda12]{Odaka2012CalabiKStability}
Yuji Odaka.
\newblock {The Calabi conjecture and K-stability}.
\newblock {\em Int. Math. Res. Not. IMRN}, (10):2272--2288, 2012.
\newblock Published version of arXiv:1010.3597.

\bibitem[Oda13]{Odaka2013}
Yuji Odaka.
\newblock {The GIT stability of Polarized Varieties via Discrepancy}.
\newblock {\em Annals of Mathematics}, 177(3):645--661, 2013.

\bibitem[Oda14]{Odaka2014}
Yuji Odaka.
\newblock {Tropical Geometric Compactification of Moduli, I – $M_g$ case}.
\newblock {\em arXiv:1406.7772}, 2014.

\bibitem[Oda15]{Odaka2015CompactModuli}
Yuji Odaka.
\newblock {Compact Moduli Spaces of K\"ahler--Einstein Fano Varieties}.
\newblock {\em Publications of the Research Institute for Mathematical Sciences}, 51(3):549--565, 2015.

\bibitem[Oda19]{Odaka2019}
Yuji Odaka.
\newblock {Tropical Geometric Compactification of Moduli, II: $A_g$ Case and Holomorphic Limits}.
\newblock {\em International Mathematics Research Notices}, 2019(21):6614--6660, 2019.

\bibitem[Oda22a]{Odaka2022b}
Yuji Odaka.
\newblock {Degenerated Calabi–Yau varieties with infinite components, moduli compactifications, and limit toroidal structures}.
\newblock {\em European Journal of Mathematics}, 3(3):1105--1157, 2022.

\bibitem[Oda22b]{Odaka2022a}
Yuji Odaka.
\newblock {PL Density Invariant for Type II Degenerating K3 Surfaces, Moduli Compactification and Hyper-Kähler Metric}.
\newblock {\em Nagoya Mathematical Journal}, 247:574--614, 2022.

\bibitem[Oda24]{Odaka2024Bubbling}
Yuji Odaka.
\newblock {Algebraic geometry of bubbling {K}\"ahler metrics}.
\newblock {\em arXiv:2406.14518}, 2024.
\newblock Preprint.

\bibitem[OO21]{OdakaOshima2021}
Yuji Odaka and Yoshiki Oshima.
\newblock {\em {Collapsing K3 Surfaces, Tropical Geometry and Moduli Compactifications of Satake, Morgan-Shalen Type}}, volume~40 of {\em Mathematical Society of Japan Memoirs}.
\newblock Mathematical Society of Japan, 2021.

\bibitem[OSS16]{odaka2016compact}
Yuji Odaka, Cristiano Spotti, and Song Sun.
\newblock {Compact moduli spaces of del Pezzo surfaces and K{\"a}hler--Einstein metrics}.
\newblock {\em Journal of Differential Geometry}, 102(1):127--172, 2016.

\bibitem[Ozu22a]{ozuch2022noncollapsed1}
Tristan Ozuch.
\newblock {Noncollapsed degeneration of Einstein 4-manifolds I}.
\newblock {\em Geometry \& Topology}, 26(4):1483--1528, 2022.

\bibitem[Ozu22b]{ozuch2022noncollapsed2}
Tristan Ozuch.
\newblock {Noncollapsed degeneration of Einstein 4-manifolds II}.
\newblock {\em Geometry \& Topology}, 26(4):1529--1634, 2022.

\bibitem[Spo14]{spotti2014deformations}
Cristiano Spotti.
\newblock {Deformations of nodal K\"ahler–Einstein Del Pezzo surfaces with discrete automorphism groups}.
\newblock {\em Journal of the London Mathematical Society}, 89(2):539--558, 2014.

\bibitem[Spo19]{spotti2019kaehler}
Cristiano Spotti.
\newblock {K\"ahler–Einstein metrics via moduli continuity method}.
\newblock In G.~Codogni, R.~Dervan, and F.~Viviani, editors, {\em Moduli of K-stable Varieties}, volume~31 of {\em Springer INdAM Series}, pages 141--152. Springer International Publishing, 2019.

\bibitem[SS17]{spotti2017explicit}
Cristiano Spotti and Song Sun.
\newblock {Explicit Gromov–Hausdorff compactifications of moduli spaces of K\"ahler–Einstein Fano manifolds}.
\newblock {\em Pure and Applied Mathematics Quarterly}, 13(3):477--515, 2017.

\bibitem[SSW20]{SongSturmWang2020}
Jian Song, Jacob Sturm, and Xiaowei Wang.
\newblock {Riemannian Geometry of K{\"a}hler--Einstein Currents III: Compactness of K{\"a}hler--Einstein Manifolds of Negative Scalar Curvature}.
\newblock {\em arXiv preprint arXiv:2003.04709}, 2020.

\bibitem[SSY16]{SpottiSunYao2016}
Cristiano Spotti, Song Sun, and Chengjian Yao.
\newblock {Existence and deformations of K\"ahler--Einstein metrics on smoothable Q-Fano varieties}.
\newblock {\em Duke Mathematical Journal}, 165(16):3043--3083, 2016.

\bibitem[Sun25]{SSBubbling}
Song Sun.
\newblock {Bubbling of {K}\"ahler-{E}instein metrics}.
\newblock {\em Pure Appl. Math. Q.}, 21(3):1317--1348, 2025.

\bibitem[Suv12]{Suvaina2012}
Ioana Suvaina.
\newblock {ALE Ricci-flat K{\"a}hler metrics and deformations of quotient surface singularities}.
\newblock {\em Annals of Global Analysis and Geometry}, 41(1):109--123, 2012.

\bibitem[SYZ96]{StromingerYauZaslow1996}
Andrew Strominger, Shing-Tung Yau, and Eric Zaslow.
\newblock {Mirror symmetry is {T}-duality}.
\newblock {\em Nuclear Physics B}, 479(1-2):243--259, 1996.

\bibitem[SZ19]{SunZhang2019}
Song Sun and Ruobing Zhang.
\newblock {Small Complex Structure Degeneration of Calabi–Yau Metrics}.
\newblock {\em arXiv preprint arXiv:1906.03368}, 2019.

\bibitem[SZ23]{sun2023no}
Song Sun and Junsheng Zhang.
\newblock {No semistability at infinity for Calabi–Yau metrics asymptotic to cones}.
\newblock {\em Inventiones Mathematicae}, 233(1):461--495, 2023.

\bibitem[SZ24]{SunZhang2024}
Song Sun and Ruobing Zhang.
\newblock {Collapsing geometry of hyperkähler 4-manifolds and applications}.
\newblock {\em Acta Mathematica}, 232(2):219--268, 2024.

\bibitem[Sze19]{S19}
Gabor Szekelyhidi.
\newblock {Degenerations of \( \mathbb{C}^n \) and Calabi–Yau metrics}.
\newblock {\em Duke Mathematical Journal}, 168(14):2651--2700, 2019.

\bibitem[Sze24]{szekelyhidi2024singular}
Gabor Szekelyhidi.
\newblock {Singular K\"ahler–Einstein metrics and RCD spaces}.
\newblock {\em arXiv preprint arXiv:2408.10747}, 2024.

\bibitem[Tam22]{Tambasco2022}
Salvatore Tambasco.
\newblock {On the continuity of Weil-Petersson volumes of the moduli space of weighted points on the projective line}.
\newblock {\em Complex Manifolds}, 9(1):206--222, 2022.

\bibitem[Thu98]{Thurston1998}
William~P. Thurston.
\newblock Shapes of polyhedra and triangulations of the sphere.
\newblock {\em Geometry \& Topology Monographs}, 1:511--549, 1998.
\newblock Part of *The Epstein Birthday Schrift*.

\bibitem[Tia90]{tian1990calabi}
Gang Tian.
\newblock {On Calabi's conjecture for complex surfaces with positive first {Chern} class}.
\newblock {\em Inventiones Mathematicae}, 101(1):101--172, 1990.

\bibitem[Tos15]{Tosatti2015}
Valentino Tosatti.
\newblock {Families of Calabi–Yau manifolds and canonical singularities}.
\newblock {\em International Mathematics Research Notices}, 2015(20):10586--10594, 2015.

\bibitem[TY90]{tian1990complete}
Gang Tian and Shing-Tung Yau.
\newblock {Complete K\"ahler manifolds with zero Ricci curvature. I}.
\newblock {\em Journal of the American Mathematical Society}, 3(3):579--609, 1990.

\bibitem[XZ24]{xu2024boundedness}
Chenyang Xu and Ziquan Zhuang.
\newblock {Boundedness of log Fano cone singularities and discreteness of local volumes}.
\newblock {\em arXiv preprint arXiv:2404.17134}, 2024.

\bibitem[Yau78]{Yau1978}
Shing-Tung Yau.
\newblock {On the Ricci curvature of a compact K{\"a}hler manifold and the complex Monge--Amp{\`e}re equation, {I}}.
\newblock {\em Communications on Pure and Applied Mathematics}, 31(3):339--411, 1978.

\bibitem[Zha15]{ZhangKE}
Yuguang Zhang.
\newblock {Collapsing of negative {K}\"ahler-{E}instein metrics}.
\newblock {\em Math. Res. Lett.}, 22(6):1843--1869, 2015.

\bibitem[Zha16]{ZhangCY}
Yuguang Zhang.
\newblock {Completion of the moduli space for polarized {C}alabi-{Y}au manifolds}.
\newblock {\em J. Differential Geom.}, 103(3):521--544, 2016.

\end{thebibliography}

\end{document}